\documentclass[12pt]{article}
\usepackage{amsmath}
\usepackage{amssymb}
\usepackage{amsthm}
\usepackage{cite}
\usepackage{epsfig}
\usepackage{graphicx}
\usepackage{subfigure}
\usepackage{psfrag}
\usepackage[latin1]{inputenc}

\newtheorem{theorem}{Theorem}[section]
\newtheorem{conjecture}[theorem]{Conjecture}

\theoremstyle{remark}
\newtheorem{remark}[theorem]{Remark}

\begin{document}

\title{\textbf{Tables, bounds and graphics of
    the smallest known sizes of complete arcs in the plane $\mathrm{PG}(2,q)$ for
    all $q\le160001$ and sporadic $q$ in the interval $[160801\ldots 430007]$}\thanks{The research
of D. Bartoli was supported by the European Community under a
Marie-Curie Intra-European Fellowship (FACE project: number
626511). The research of  G. Faina, S. Marcugini, and F.
Pambianco was
 supported in part by Ministry for Education, University
and Research of Italy (MIUR) (Project ``Geometrie di Galois e
strutture di incidenza'')
 and by the Italian National Group for Algebraic and Geometric Structures and their Applications
 (GNSAGA - INDAM).
 The research of A.A.~Davydov and A.A.~Kreshchuk was carried out at the IITP RAS at the expense of the Russian
Foundation for Sciences (project 14-50-00150).}}
\date{}
\author{Daniele Bartoli \\
{\footnotesize Department of  Mathematics, Ghent University,}\\
{\footnotesize Krijgslaan 281-B, Gent 9000, Belgium.
E-mail: dbartoli@cage.be}\\
 \and Alexander A.
Davydov \\
{\footnotesize Institute for Information Transmission Problems
(Kharkevich
institute), Russian Academy of Sciences}\\
{\footnotesize Bol'shoi Karetnyi per. 19, GSP-4, Moscow,
127994, Russia. E-mail: adav@iitp.ru}
\and Giorgio Faina \\
{\footnotesize Dipartimento di Matematica e Informatica,
Universit\`{a}
degli Studi di Perugia, }\\
{\footnotesize Via Vanvitelli~1, Perugia, 06123, Italy. E-mail:
 giorgio.faina@unipg.it} \and Alexey A.
Kreshchuk {\footnotesize \ } \\
{\footnotesize Institute for Information Transmission Problems
(Kharkevich
institute), Russian Academy of Sciences}\\
{\footnotesize Bol'shoi Karetnyi per. 19, GSP-4, Moscow,
127994, Russia. E-mail: krsch@iitp.ru}
\and Stefano Marcugini and Fernanda Pambianco \\
{\footnotesize Dipartimento di Matematica e Informatica,
Universit\`{a}
degli Studi di Perugia, }\\
{\footnotesize Via Vanvitelli~1, Perugia, 06123, Italy. E-mail:
\{stefano.marcugini,fernanda.pambianco\}@unipg.it}}
\maketitle

\noindent\textbf{Abstract.} In the projective planes
$\mathrm{PG}(2,q)$, we collect \textbf{the smallest known sizes
of complete arcs} for the regions
\begin{align*}
&\mbox{ \textbf{all} } q\le160001,~~ q \mbox{ prime power};\\
&Q_{4}=\{34 \mbox{ sporadic }q'\mbox{s in the interval }[160801\ldots430007]
, \mbox{ see Table 3}\}.
\end{align*}
 For
$q\le160001$, the \textbf{collection of arc sizes is complete}
in the sense that arcs for \textbf{all prime powers} are
considered. This proves new upper bounds on the smallest size
$t_{2}(2,q)$ of a complete arc in $\mathrm{PG}(2,q)$, in
particular
    \begin{align}
t_{2}(2,q)&<0.998\sqrt{3q\ln q}<1.729\sqrt{q\ln q}&\mbox{ for }&&7&\le q\le160001;\label{eqAbs_1} \\
t_{2}(2,q)&<\sqrt{q}\ln^{0.7295}q&\mbox{ for }&&109&\le q\le160001;\label{eqAbs_2}\\
t_{2}(2,q)&<\sqrt{q}\ln^{c_{up}(q)}q,~~
c_{up}(q)=\frac{0.27}{\ln q}+0.7,&\mbox{ for }&&19&\le q\le160001;\label{eqAbs_3}\\
t_{2}(2,q)&<0.6\sqrt{q}\ln^{\varphi_{up}(q;0.6)} q,~~
\varphi_{up}(q;0.6)=\frac{1.5}{\ln q}+0.802,&\mbox{ for }&&19&\le q\le160001.\label{eqAbs_4}
\end{align}
Moreover, the bounds \eqref{eqAbs_2} -- \eqref{eqAbs_4} hold
also for $q\in Q_{4}$. Also,
    \begin{align}
t_{2}(2,q)&<1.006\sqrt{3q\ln q}<1.743\sqrt{q\ln q}&\mbox{ for }&&q\in Q_{4}.\label{eqAbs_5}
\end{align}
Our investigations and results allow to conjecture that the
bounds \eqref{eqAbs_2} -- \eqref{eqAbs_5} hold for all
$q\geq109$.



\numberwithin{equation}{section}
\section{Introduction. The main results}\label{sec_Intro}

Let $\mathrm{PG}(2,q)$ be the projective plane over the Galois
field $\mathbb{F}_{q}$ of $q$ elements. An $n$-arc is a set of
$n$ points no three of which are collinear. An $n$-arc is
called complete if it is not contained in an $(n+1)$-arc of
$\mathrm{PG}(2,q)$. For an introduction to projective
geometries over finite fields see \cite
{HirsBook,SegreLeGeom,SegreIntrodGalGeom}.

 The relationship among the
theory of $n$-arcs, coding theory, and mathematical statistics
is presented in \cite{HirsSt-old,HirsStor-2001}, see also
\cite{Landjev,LandjevStorme2012}. In particular, a complete arc
in a plane $\mathrm{PG}(2,q),$ the points of which are treated
as 3-dimensional $q$-ary columns, defines a parity check matrix
of a $q$-ary linear code with codimension 3, Hamming distance
4, and covering radius 2. Arcs can be interpreted as linear
maximum distance separable (MDS) codes \cite[Sec.~7]{szoT93},
\cite{thaJ92d} and they are related to optimal coverings arrays
\cite{Hartman-Haskin}, superregular matrices \cite {Keri}, and
quantum codes \cite{BMP-quantum}.

A point set $S\subset \mathrm{PG}(2,q)$ is 1-\emph{saturating}
if any point of $ \mathrm{PG}(2,q)\setminus S$ is collinear
with two points in $S$ \cite
{BFMP-JG2013,DGMP-AMC,DMP-JCTA2003,Giul2013Survey,Ughi-sat}.
The points of a 1-saturating set in $\mathrm{PG}(2,q)$ form a
parity check matrix of a \emph{linear covering code} with
codimension 3,  Hamming distance 3 or 4, and covering radius 2.
An open problem is to find small 1-saturating sets
(respectively, short covering codes). A complete arc in
$\mathrm{PG}(2,q)$ is, in particular, a 1-saturating set; often
the smallest known complete arc is the smallest known
1-saturating set \cite
{BFMP-JG2013,DGMP-AMC,DMP-JCTA2003,Giul2013Survey,Pace-A5}. Let
$\ell_{1}(2,q)$ be the smallest size of a 1-saturating set in
$\mathrm{PG}(2,q)$. In \cite {BorSzonTicDefinSets,Kovacs}, for
$q$
 large enough, by probabilistic methods the following
upper bound  is proved (for $3\sqrt{2}$ see~\cite[p.\thinspace
24]{BorSzonTicDefinSets}):
\begin{equation}  \label{eq1_1sat}
\ell_{1}(2,q)<3\sqrt{2}\sqrt{q\ln q}<5\sqrt{q\ln q}.
\end{equation}

Let $t_{2}(2,q)$ be the \textbf{smallest size of a complete arc
in $\mathrm{PG}(2,q)$}.

One of the main problems in the study of projective planes,
which is also of interest in coding theory, is the finding of
the spectrum of possible sizes of complete arcs. In particular,
the value of $t_{2}(2,q)$ is interesting. Finding estimates of
the minimum size $t_{2}(2,q)$ is a hard open problem, see e.g.
\cite[Sec.\,4.10]{HirsThas-2015}.

This paper is devoted to \emph{upper bounds} for $t_{2}(2,q)$.

Surveys of results on the sizes of plane complete arcs, methods
of their construction, and the comprehension of the relating
properties can be found in\cite%
{BDFKMP-PIT2014,BDFKMP-JGtoappear,BDFMP-DM,BDFMP-JG2013,BDFMP-JG2015,BDMP-3cycle,%
DFMP-JG2005,DFMP-JG2009,FP,Giul2013Survey,%
HirsSurvey83,HirsBook,HirsSt-old,HirsStor-2001,HirsThas-2015,KV,LombRad,pelG77,pel93,%
SegreLeGeom,SegreIntrodGalGeom,SZ,szoT87a,szoT89survey,szoT93,Szonyi1997surveyCurves}.

Problems connected with small complete arcs in $\mathrm{PG}(2,q)$ are considered in \cite%
{abaV83,AliPhD,Ball-SmallArcs,BDFKMP-Bulg2013,BDFKMP-ArXiv2013,BDFKMP-ArXiv,BDFKMP-ACCT2014Conject,%
BDFKMP-JGtoappear,BDFMP-DM,BDFMP-JG2013,ComputBound-Svetlog2014,BDFKMP-PIT2014,%
BDFKMP-ArXiv2015lexi,BDFMP-Bulg2012a,BDFMP-JG2015,BDMP-Bulg2012b,BFMP-JG2013,%
BFMPD-ENDM2013,BDFMP-ArXivFOP,BDFMP-ArXivRandom,Blokhuis,CoolStic2009,CoolStic2011,%
DFMP-JG2005,DFMP-JG2009,DGMP-Innov,DGMP-JCD,%
DGMP-AMC,DMP-JCTA2003,FainaGiul,FMMP-1977,FP,FPDM,GacsSzonyi,Giul2000,%
Giul2007affin,Giul2007even,Giul2013Survey,GKMP-A6invar,GiulUghi,Gordon,Hadnagy,HirsSurvey83,%
HirsBook,HirsSad,HirsSt-old,HirsStor-2001,HirsThas-2015,KV,korG83a,%
LisPhD,LisMarcPamb2008,LombRad,MMP-q29,MMP-q25,Ost,Pace-A5,pelG77,%
pel93,Polv,SegreLeGeom,SegreIntrodGalGeom,SZ,szoT87a,szoT87b,%
szoT89survey,szoT93,Szonyi1997surveyCurves,Ughi-sqrt-log2,UghiAlmost,Voloch87,Voloch90}%
, see also references therein.

The exact values of $t_{2}(2,q)$ are known only for $q\leq 32$;
see \cite
{AliPhD,CoolStic2009,CoolStic2011,FMMP-1977,Gordon,HirsBook,HirsSad,%
HirsSt-old,HirsStor-2001,HirsThas-2015,MMP-q29,MMP-q25}
and the paper \cite {BFMP-JG2013} where the equalities
$t_{2}(2,31)=t_{2}(2,32)=14$ are established.

The following lower bounds hold (see \cite
{Ball-SmallArcs,Blokhuis,Polv,SegreLeGeom} and references
therein):
\begin{equation}
t_{2}(2,q)>\left\{
\begin{array}{ll}
\sqrt{2q}+1 & \text{for any }q, \\
\sqrt{3q}+\frac{1}{2} & \text{for }q=p^{h},\text{ }p\text{ prime, }h=1,2,3.
\end{array}
\right.   \label{eq1_trivlower}
\end{equation}

Let $t(\mathcal{P}_{q})$ be the size of the smallest complete
arc in any (not necessarily Desarguesian) projective plane
$\mathcal{P}_{q}$ of order $q$. In \cite{KV}, for $q$
\emph{large enough}, the following result is proved by
\emph{probabilistic methods} (we take into account the remark
\cite[p.\ 320]{KV} that all logarithms in \cite{KV} have the
natural base):
\begin{equation}
t(\mathcal{P}_{q})\leq \sqrt{q}\ln^{C}q,\text{ }C\leq 300,
\label{eq1_KimVu_c=300}
\end{equation}
where $C$ is a constant independent of $q$ (so-called universal
or absolute constant). The authors of  \cite{KV} conjecture
that the constant $C$ can be reduced to $C=10$. A survey and
analysis of random constructions for geometrical objects can be
found in \cite {BorSzonTicDefinSets,GacsSzonyi,Kovacs}.

In $\mathrm{PG}(2,q)$, complete arcs are obtained by
\emph{algebraic constructions} (see
\cite[p.\thinspace209]{HirsStor-2001}) with sizes approximately
$\frac{1}{3}q$ \cite
{abaV83,BDFMP-DM,korG83a,SZ,szoT87a,szoT87b,Voloch90},
$\frac{1}{4}q$ \cite {BDFMP-DM,korG83a,szoT89survey},
$2q^{0.9}$ for $q>7^{10}$ \cite{SZ}, and $ 2.46q^{0.75} \ln q$
for big prime $q$ \cite{Hadnagy}. It is noted in \cite[
Sec.~8]{GacsSzonyi} that the smallest size of a complete arc in
$ \mathrm{PG}(2,q)$ obtained via algebraic constructions is
\begin{equation}\label{eq1_cq34}
cq^{3/4}
\end{equation}
 where $c$ is a universal constant  \cite[ Sec.\ 3]{szoT89survey}, \cite[ Th.\
6.8]{szoT93}.

Thus, there is a substantial gap between the known upper bounds
and the lower bounds on $t_{2}(2,q)$, see
\eqref{eq1_trivlower}--\eqref{eq1_cq34}. The gap is essentially
reduced if one consider the lower bound \eqref{eq1_trivlower}
for complete arcs and the upper bound \eqref{eq1_1sat} for
1-saturating sets. However, though complete arcs are
1-saturating sets they represent a narrower class of objects.
Therefore, for complete arcs, one may not use the bound
\eqref{eq1_1sat} directly. Nevertheless, the common nature of
complete arcs and 1-saturating sets allows to hope for upper
bounds on $t_{2}(2,q)$ similar to \eqref{eq1_1sat}. The hope is
supported by numerous experimental data, see below.

 In  \cite[p.\ 313]{KV}, it is noted (with reference to the work \cite{Blokhuis}) that in a
preliminary report of 1989, J.C. Fisher  obtained by computer
search complete arcs in many planes of small orders and
conjectured that the average size of a complete arc in
$\mathrm{PG}(2,q)$ is about $\sqrt{ 3q\log q}$.

In \cite{BDFKMP-PIT2014}, see also
\cite{BDFKMP-ACCT2014Conject}, an attempt to obtain a
theoretical upper bound on $t_{2}(2,q)$ with the main term of
the form $c\sqrt{q\ln q}$, where $c$ is a small universal
constant, is done. The reasonings of \cite{BDFKMP-PIT2014} are
based on explanation of the working mechanism of a step-by-step
greedy algorithm for constructing complete arcs in
$\mathrm{PG}(2,q)$ and on quantitative estimations of the
algorithm. For more than half of the steps of the iterative
process, these estimations are proved rigorously. The natural
(and well-founded) conjecture that they hold for the rest of
steps is done, see \cite[Conject.\ 2]{BDFKMP-PIT2014}. As a
result the following conjectural upper bounds are formulated:

\begin{conjecture}\emph{\textbf{\cite{BDFKMP-PIT2014}}}\label{conj1}
In $\mathrm{PG}(2,q)$, under conjecture given in
\emph{\cite[Conject.\ 2]{BDFKMP-PIT2014}}, the following
estimates hold:
    \begin{align}
&t_{2}(2,q)<\sqrt{q}\sqrt{3\ln q+\ln \ln q+\ln
3}+\sqrt{\frac{q}{3\ln q}}+3,
 \label{eq1_Prob bounds}\displaybreak[3]\\
  &t_{2}(2,q)<1.87\sqrt{q\ln q}<1.08\sqrt{3q\ln q}.\label{eq1_!<}
  \end{align}
\end{conjecture}
Moreover, in \cite{BDFKMP-PIT2014} it is conjectured that the
upper bounds \eqref{eq1_Prob bounds}, \eqref{eq1_!<} hold for
all $q$ without any extra conditions.

 Denote by $\overline{t}_{2}(2,q)$ \textbf{the smallest
\emph{known }size of a complete arc in} $\mathrm{PG}(2,q)$.
Clearly,
\begin{equation}\label{eq1_main-ineq}
t_{2}(2,q)\le\overline{t}_{2}(2,q).
\end{equation}

From  \cite{BDFMP-DM,BDFMP-JG2013}, see also references
therein, it follows that  $\overline{ t}_{2}(2,q)<4\sqrt{q}$
$\mbox{ for }$ $q\leq 841$ and $q=857$. Complete $(4
\sqrt{q}-4)$-arcs are obtained for odd $q=p^{2}$ with $p\leq
41$ and $p=7^{2}$ \cite {Giul2000,GiulUghi} and for
$q=2^{6},2^{8},2^{10}$ \cite{DGMP-JCD,FainaGiul}. So,
\begin{equation}
t_{2}(2,q)<4\sqrt{q}\quad \text{for }q\leq 841,\text{ }
q=857,31^{2},2^{10},37^{2},41^{2},7^{4}.\notag
\end{equation}
For $q\leq 151,$ a number of improvements of
$\overline{t}_{2}(2,q)$, in comparison with
\cite{BDFMP-DM,BDFMP-JG2013}, are given in \cite{Pace-A5} and
cited in \cite{BDFKMP-ArXiv}.

For $q\leq 13627$ and sporadic $q\le 45893$, the values of
$\overline{t} _{2}(2,q)$ (up to January 2013) are collected in
 \cite {BDFMP-DM,BDFMP-JG2013} where the following
results are obtained:
\begin{align}
& t_{2}(2,q)<4.5\sqrt{q} \text{ for }q\leq 2647,\text{ }
q=2659,2663,2683,2693,2753,2801; \displaybreak[3]   \notag\\
& t_{2}(2,q)<5\sqrt{q} \text{ for }q\leq 9497,\text{ }
q=9539,9587,9613,9623,9649,9689,\displaybreak[3]\notag \\
&\phantom{t_{2}(2,q)<5\sqrt{q} \text{ for }q\leq 9497,\text{ }
q=\,\,}9923,9973. \notag
\end{align}

Let $Q_{1}$, $Q_{2}$, $Q_{3}$, $Q_{4}$, and $Q$ be the
following sets of values of $q$:
\begin{align}
Q_{1}&=\{2\le q\le49727,~ q\mbox{ prime power}\}\cup\{2\le q\leq
150001,~ q \mbox{ prime}\}\,\cup  \displaybreak[0]\notag\\*
& \{40 \mbox{ sporadic prime $q$'s in the interval
}[150503\ldots410009]\};
 \label{eq1_Q1}\displaybreak[3]\\
   \label{eq1_Q2}
Q_{2}&=\{49727< q\le 150001,~q=p^{h},~h\ge2,~p \mbox{ prime}\}\,\cup \{
430007\};\displaybreak[3]\\
Q_{3}&=\{150001< q\le 160001,~q\mbox{ prime power}\};
   \label{eq1_Q3}\displaybreak[3]\\
Q_{4}&=\{34 \mbox{ sporadic }q'\mbox{s in the interval }[160801\ldots430007], \mbox{ see Table 3}\};
   \label{eq1_Q4}\displaybreak[3]\\
Q&=Q_{1}\cup Q_{2}\cup Q_{3}\cup Q_{4}=\{2\le q\le160001, q \mbox{ prime power}\}\,\cup\notag\\*
&\{34 \mbox{ sporadic }q'\mbox{s in the interval }[160801\ldots430007]\}.\label{eq1_Q}
\end{align}
For $q\in Q_{1}$, the
 values of $\overline{t}_{2}(2,q)$ (up to August
2014) are collected in the work \cite{BDFKMP-ArXiv}. The values
of $q$ in the interval $[150503\ldots410009]$, corresponding to
the 2-nd row of the formula~\eqref{eq1_Q1}, are given in
\cite[Tab.\,6]{BDFKMP-ArXiv}. For $q\in Q_{2}$, the values  of
$\overline{t}_{2}(2,q)$ are obtained in \cite{BDFKMP-PIT2014}
and collected in \cite[Tab.\,1]{BDFKMP-JGtoappear} and Table 2
of this work. For  $q\in Q_{3}$ the values of
$\overline{t}_{2}(2,q)$ are obtained in
\cite{BDFKMP-JGtoappear} and collected in Table 6 of this work.
The values of $q$ in the interval $[160801\ldots430007]$
corresponding to $Q_{4}$ and the 2-nd row of the
formula~\eqref{eq1_Q} are given in Table~3 of this work. They
are obtained in \cite{BDFKMP-ArXiv,BDFKMP-JGtoappear} and in
this work.\footnote{In the works
\cite{BDFKMP-ArXiv2013,BDFKMP-ArXiv,BDFKMP-PIT2014,BDFKMP-JGtoappear}
and in the present work, the most of calculations are done
using computational resources of Multipurpose Computing Complex
of National Research Centre ``Kurchatov Institute'',
http://computing.kiae.ru} For $q=2^{18}$ the result of
\cite{DGMP-Innov} is used, see~\eqref{eq1_even} below. Note
that $160801=401^{2}$.

From results collected in \cite{BDFKMP-ArXiv}, we have
\begin{align}
t_{2}(2,q)<5.5\sqrt{q} \text{ for }q\leq 38557,\text{ }
q\neq 36481,37537,37963,38039,38197.  \notag
\end{align}

For even $q=2^{h}$, small complete arcs in planes are a base
for inductive infinite families of small complete caps in the
projective spaces $\mathrm{PG}(N,q)$, see \cite{DGMP-JCD}. For
 $h\leq 15$, the smallest known sizes of complete arcs in
$\mathrm{PG}(2,2^{h})$ are collected in \cite{BDFKMP-ArXiv},
see also
\cite{BDFMP-DM,DFMP-JG2009,DGMP-JCD,BDFMP-JG2013,FainaGiul} and
references therein. For $q=2^{16},2^{17}$ complete arcs are
obtained in \cite{BDFKMP-PIT2014}. Finally, complete
${(6\sqrt{q}-6)}${-arcs} in $\mathrm{PG}(2,4^{2h+1})$, are
constructed in \cite{DGMP-Innov}; for $h\leq 4$ it is proved
that they are complete. This gives a complete $3066$-arc in
$\mathrm{PG}(2,2^{18})$. Sizes of the arcs for $h\ge11$ are as
follows:
\begin{align}
& \overline{t}_{2}(2,2^{11})=199,~\overline{t}_{2}(2,2^{12})=300,~
\overline{t}_{2}(2,2^{13})=449,~\overline{t} _{2}(2,2^{14})=665,\displaybreak[3]\notag\\
&\overline{t}_{2}(2,2^{15})=987,~\overline{t}_{2}(2,2^{16})=1453,~
\overline{t}_{2}(2,2^{17})=2141,~\overline{t}_{2}(2,2^{18})=3066.\label{eq1_even}
\end{align}

Note that, for $q\in Q$, in the works
\cite{BDFKMP-ArXiv2013,BDFKMP-ArXiv,BDFKMP-PIT2014,%
BDFKMP-JGtoappear,BDFMP-DM,BDFMP-JG2013} and in the present
paper, most of the values $\overline{t }_{2}(2,q)$ have been
obtained by computer search using  randomized greedy algorithms
described in \cite
{BDFKMP-JGtoappear,BDFMP-DM,BDFMP-JG2013,BDFMP-JG2015,DFMP-JG2005,DFMP-JG2009,DMP-JG2004}.
The step-by-step greedy algorithm adds to
    the arc in every step a point providing the maximum
    possible (for the given step) number of new covered
    points.

Another way for obtaining of small plane complete arcs is
applying a step-by-step algorithm with fixed order of points
(FOP), see
\cite{BDFKMP-ArXiv2015lexi,ComputBound-Svetlog2014,BDFKMP-JGtoappear,BDFMP-JG2015,%
BDMP-Bulg2012b,BFMPD-ENDM2013,BDFMP-ArXivFOP}. The algorithm
FOP fixes a particular order on points of $\mathrm{PG}(2,q)$.
In every step the algorithm FOP adds to an incomplete running
arc the next point in this order  not lying on bisecants of
this arc. In this work, we do not consider the results obtained
by the FOP algorithm.

It should be emphasized that in the present work we collect and
use the computer search results for \textbf{all prime power}
$q$ in the region $q\le 160001$. In this sense, we say that the
\textbf{computer search} in this region and the \textbf{arc
sizes collection} are \textbf{complete}.

 In the works of the authors
\cite{BDFKMP-Bulg2013,BDFMP-JG2013,BDFMP-JG2015,BFMPD-ENDM2013,BDMP-Bulg2012b},
non-standard types of upper bounds on $t_{2}(2,q)$ are
proposed. They are based on \emph{decreasing functions} $c(q)$
and
 $\varphi(q;D)$ defined by the following relations:
\begin{align}\label{eq1 c(q)}
&t_{2}(2,q)=\sqrt{q}\ln ^{c(q)} q,~~c(q) = \frac{\ln(t_{2}(2,q)/\sqrt{q})}{\ln \ln q},\displaybreak[3]\\
&t_{2}(2,q)=D\sqrt{q} \ln^{\varphi(q;D)} q,~~\varphi(q;D) = \frac{\ln(t_{2}(2,q)/D\sqrt{q})}{\ln \ln q},\label{eq1 phi(q)}
\end{align}
where $D$ is a universal constant independent of $q$.

Also, in the paper \cite{BDFKMP-JGtoappear}, the function
$h(q)$ is defined so that
\begin{align}
&t_{2}(2,q)=h(q)\sqrt{3q \ln q}.\label{eq1_h(q)}
\end{align}

The \emph{non-standard types of upper bounds} estimate the
functions $c(q)$, $\varphi(q;D)$, and $h(q)$ instead of
$t_{2}(2,q)$. This allows to do estimates and graphics more
expressive, see Figures\linebreak \ref{fig1_t_Gr&bounds} --
\ref{fig11_varphi_Gr&bounds} below.

We denote the following set of values of $q$:
\begin{align}
Q_{4}^{*}=\{q\,:\,q\in Q_{4}, q\le190027\}\cup\{q=2^{18},380041\}\subset Q_{4}.\label{eq1_Q4star}
\end{align}

The following theorem summarizes the main results of this paper
using the mentioned types of upper bounds. The theorem is based
on the complete computer search the results of which are
collected in
\cite{BDFKMP-ArXiv2013,BDFKMP-ArXiv,BDFMP-DM,BDFMP-JG2013} and
in Tables 1 -- 6 of the present paper.

\begin{theorem}\label{th1_main}
Let $t_{2}(2,q)$ be the smallest size of a complete arc in the
projective plane $\mathrm{PG}(2,q)$. Let $Q_{4}$ and
$Q_{4}^{*}$ be the sets of values of $q$ given by relations
\eqref{eq1_Q4}, \eqref{eq1_Q4star} and Table \emph{3}. In
$\mathrm{PG}(2,q)$ the following hold.
\begin{itemize}
   \item[A.] There are the following upper bounds with a
       \textbf{fixed degree of logarithm} of~$q$\linebreak
       (\textbf{FDL-bounds}):
\begin{itemize}
\item[(i)]
    \begin{align}
&t_{2}(2,q)<
0.998\sqrt{3q\ln q}<1.729\sqrt{q\ln q}~~\mbox{ for }~~7\le q\le160001\mbox{ and }q\in Q_{4}^{*},\label{eq1_Bounds1G}\displaybreak[0]\\
&t_{2}(2,q)<
1.006\sqrt{3q\ln q}<1.743\sqrt{q\ln q}~~\mbox{ for }~~7\le q\le160001\mbox{ and }q\in Q_{4}. \label{eq1_Bounds1GQ4}
\end{align}
\item[(ii)]
    \begin{align}
&t_{2}(2,q)<\sqrt{q}\ln^{0.7295}q~~\mbox{ for }~~109\le q\le160001\mbox{ and }q\in Q_{4}.\label{eq1_Bounds2G}
\end{align}
        \end{itemize}
        \item[B.] There are the following upper bounds with
            a \textbf{decreasing degree of logarithm} of
            $q$ (\textbf{DDL-bounds}):
\begin{itemize}
\item[(i)] Let $t_{2}(2,q)=\sqrt{q}\ln^{c(q)}q$. We
    have
    \begin{align}
&c(q)<\frac{0.27}{\ln q}+0.7~~\mbox{ for }~~19\le q\le160001\mbox{ and }q\in Q_{4}.\label{eq1_Bounds3G}
\end{align}
\item[(ii)] Let
    $t_{2}(2,q)=0.6\sqrt{q}\ln^{\varphi(q;0.6)} q$. We
    have
    \begin{align}
\varphi(q;0.6)<\frac{1.5}{\ln q}+0.802~~\mbox{ for }~~19\le q\le160001\mbox{ and }q\in Q_{4}.   \label{eq1_Bounds4G}
\end{align}
        \end{itemize}
        \end{itemize}
Complete arcs in $\mathrm{PG}(2,q)$ satisfying the upper bounds
\eqref{eq1_Bounds1G} -- \eqref{eq1_Bounds4G} can be constructed
with the
    help of the step-by-step greedy algorithm which adds to
    the arc in every step a point providing the maximum
    possible (for the given step) number of new covered
    points.
\end{theorem}

Calculations executed for sporadic $q\le430007$ strengthen the
confidence for validity of the bounds of Theorem
\ref{th1_main}. Also, it is important that the bounds
\eqref{eq1_Bounds1G}--\eqref{eq1_Bounds4G} are close to the
conjectural (but well-founded) bounds of \cite{BDFKMP-PIT2014},
see \eqref{eq1_Prob bounds}, \eqref{eq1_!<} in Conjecture
\ref{conj1}. On the whole, our investigations and results (see
also figures below) allow to conjecture the following.

\begin{conjecture}\label{conj1_for all q}
The upper bounds \eqref{eq1_Bounds1GQ4} -- \eqref{eq1_Bounds4G}
hold for all $q\ge109$.
\end{conjecture}

\section{Greedy algorithms. The smallest known complete arcs}
\label{sec-greedy}
\subsection{Randomized greedy algorithms}
The most of the smallest known complete arcs \emph{for
 $q\in Q$ used in the present paper} are obtained by computer search based on
randomized greedy algorithms
\cite{BDFKMP-ArXiv,BDFMP-DM,BDFMP-JG2013,BDFMP-JG2015,DFMP-JG2005,DFMP-JG2009,DMP-JG2004},
with the exception of (relatively few) theoretically
constructed complete arcs, see, for example,
\cite{DGMP-Innov,DGMP-JCD,Giul2000,GKMP-A6invar,GiulUghi,Pace-A5,FainaGiul}
and references therein.

In every step a randomized greedy  algorithm maximizes an
objective function $f$ but some steps are executed in a random
manner. The number of these steps, their ordinal numbers, and
some other parameters of the algorithm have been taken
intuitively. Also, if the same maximum of $f$ can be obtained
in distinct ways, one way is chosen randomly.

We begin to construct a complete arc by using a starting point
set $ S_{0} $. In the $i$-th step one point is added to the set
$S_{i-1}$ and we obtain a point set $S_{i}$. As the value of
the objective function $f$ we consider \emph{the number of
covered }points in $\mathrm{PG}(2,q),$ that is, points that lie
on bisecants of $S_{i}$.

On every \textquotedblleft random\textquotedblright\ $i$-th
step we take $ d_{q,i}$ \emph{randomly chosen points} of
$\mathrm{PG}(2,q)$ \emph{not covered by }$ S_{i-1}$ and compute
the objective function $f$ adding each of these $ d_{q,i} $
points to $S_{i-1}$. The point providing the maximum of $f$ is
included into~$S_{i}.$ On every \textquotedblleft
non-random\textquotedblright\ $j$-th step we consider \emph{all
points not covered by }$S_{j-1}$ and add to $S_{j-1}$ the point
providing the maximum of $f.$

 As
$S_{0}$ we can use a subset of points of an arc obtained in
previous stages of the search.

A generator of random numbers is used for a random choice. To
get arcs with distinct sizes, the starting conditions of the
generator are changed for the same set $S_{0}$. In this way the
algorithm works in a convenient limited region of the search
space to obtain examples improving the size of the arc from
which the fixed points have been taken.

In order to obtain arcs with new sizes, sufficiently many
attempts should be made with  randomized greedy algorithms.
``Predicted'' sizes  could be useful for understanding if a
good result has been obtained. For $\mathrm{PG}(2,q)$, the
predicted sizes  can be obtained in distinct ways including
upper bounds considered in Introduction and complete arc sizes
for smaller $q$ calculated before. See, for example, approaches
used in \cite{BDFMP-JG2013}. If the result is not close to the
predicted size, the attempts are continued.

We obtain small complete arcs in $\mathrm{PG}(2,q)$ in two
stages.

At the 1-st stage, we take the frame as $S_{0}$ and create a
starting complete arc $K_{0}$ using in the beginning of the
process $\delta _{q}$ \textquotedblleft
random\textquotedblright\ steps with distinct $d_{q,i}.$ All
the next steps  are \textquotedblleft
non-random\textquotedblright .

At the 2-nd stage we execute $n_{q}$ attempts to get a complete
arc. For every attempt, the starting conditions of the random
generator are other than in the previous ones. But the set
$S_{0}$ is the same. Two or three of the first five steps of
every attempt are \textquotedblleft random\textquotedblright ,
the rest of them are \textquotedblleft
non-random\textquotedblright .

The values of $d_{q,i}$, $\delta _{q}$, and $n_{q}$  are given
intuitively depending on $q$ and (for $d_{q,i})$ on $|S_{i-1}|$
and on the stage of the process. Of course, CPU performance
affects the algorithm parameters choice.

Thus, arc sizes obtained by the randomized greedy algorithms
depend on many factors. Nevertheless, we repeat that, for
 $q\in Q$, the most of the smallest
known complete arcs are obtained with the help of these
algorithms.

\subsection{Upper bounds on $t_{2}(2,q)$, $q\in Q$, based on the smallest known complete arcs}

\emph{\textbf{In this work we collected the smallest known
sizes $\overline{t}_{2}(2,q)$ of complete arcs in
$\mathrm{PG}(2,q)$ (up to June 2015) for $q\in Q$}}, see
\eqref{eq1_Q}.

 The
values of $\overline{t}_{2}(2,q)$ are written in Tables 1 -- 6,
see Appendix. We take into account the results  obtained and
collected in \cite
{abaV83,AliPhD,Ball-SmallArcs,BDFKMP-Bulg2013,BDFKMP-ArXiv2013,BDFKMP-ArXiv,%
BDFKMP-ACCT2014Conject,BDFKMP-PIT2014,BDFKMP-JGtoappear,ComputBound-Svetlog2014,%
BDFMP-Bulg2012a,BDFMP-DM,BDFMP-JG2013,%
BDFMP-JG2015,BDMP-Bulg2012b,BFMP-JG2013,BFMPD-ENDM2013,DFMP-JG2005,%
CoolStic2009,CoolStic2011,DFMP-JG2009,DGMP-Innov,DGMP-JCD,DGMP-AMC,%
FainaGiul,FMMP-1977,FP,FPDM,GacsSzonyi,%
Giul2000,Giul2007affin,Giul2007even,Giul2013Survey,GKMP-A6invar,GiulUghi,%
Gordon,Hadnagy,HirsBook,HirsSad,HirsSt-old,HirsStor-2001,HirsThas-2015,%
KV,korG83a,Landjev,LisPhD,LisMarcPamb2008,LombRad,%
MMP-q29,MMP-q25,Ost,Pace-A5,pelG77,pel93,Polv,SegreLeGeom,%
SegreIntrodGalGeom,SZ,szoT87a,szoT87b,szoT89survey,szoT93,Szonyi1997surveyCurves,%
Ughi-sqrt-log2,UghiAlmost,Voloch87,Voloch90}.

For $2\le q\le 13627$ we use arc's sizes collected in the paper
\cite{BDFMP-JG2013}. Some sizes of \cite{BDFMP-JG2013} are
improved; they are collected in Table~1. For $q\le151$, the
improvements are taken from \cite{Pace-A5}; for $q\ge2383$, the
improved sizes are obtained in \cite{BDFKMP-ArXiv2013}.

In Table 2, for $49729\leq q\leq 160801$, $q$ non-prime, the
smallest known sizes
 $\overline{t}_{2}(2,q)$ (for short $\overline{t}_{2}$)
 of complete arcs in $\mathrm{PG}(2,q)$ are collected.
 All $q\in Q_{2}$, see \eqref{eq1_Q2}, are included to this table.

In Table 3, for 34 sporadic $q$'s in the interval
$[160801\ldots
 430007]$ and for $q=124477,\linebreak 128683,160001$, the smallest known sizes
 $\overline{t}_{2}=\overline{t}_{2}(2,q)$
 of complete arcs in $\mathrm{PG}(2,q)$ are written.
 All $q\in Q_{4}$, see \eqref{eq1_Q4}, are included to this table.
 The sizes  that are new or smaller in comparison with \cite{BDFKMP-ArXiv}
 are written either in bold font (if they are obtained in \cite{BDFKMP-JGtoappear}) or in italic-bold
 one (if they are obtained in this work).
 The new results for $\overline{t}_{2}(2,124477)$ and $\overline{t}_{2}(2,128683)$ satisfy the
 bound~\eqref{eq1_Bounds1G}.

\begin{figure}[htbp]
\includegraphics[width=\textwidth]{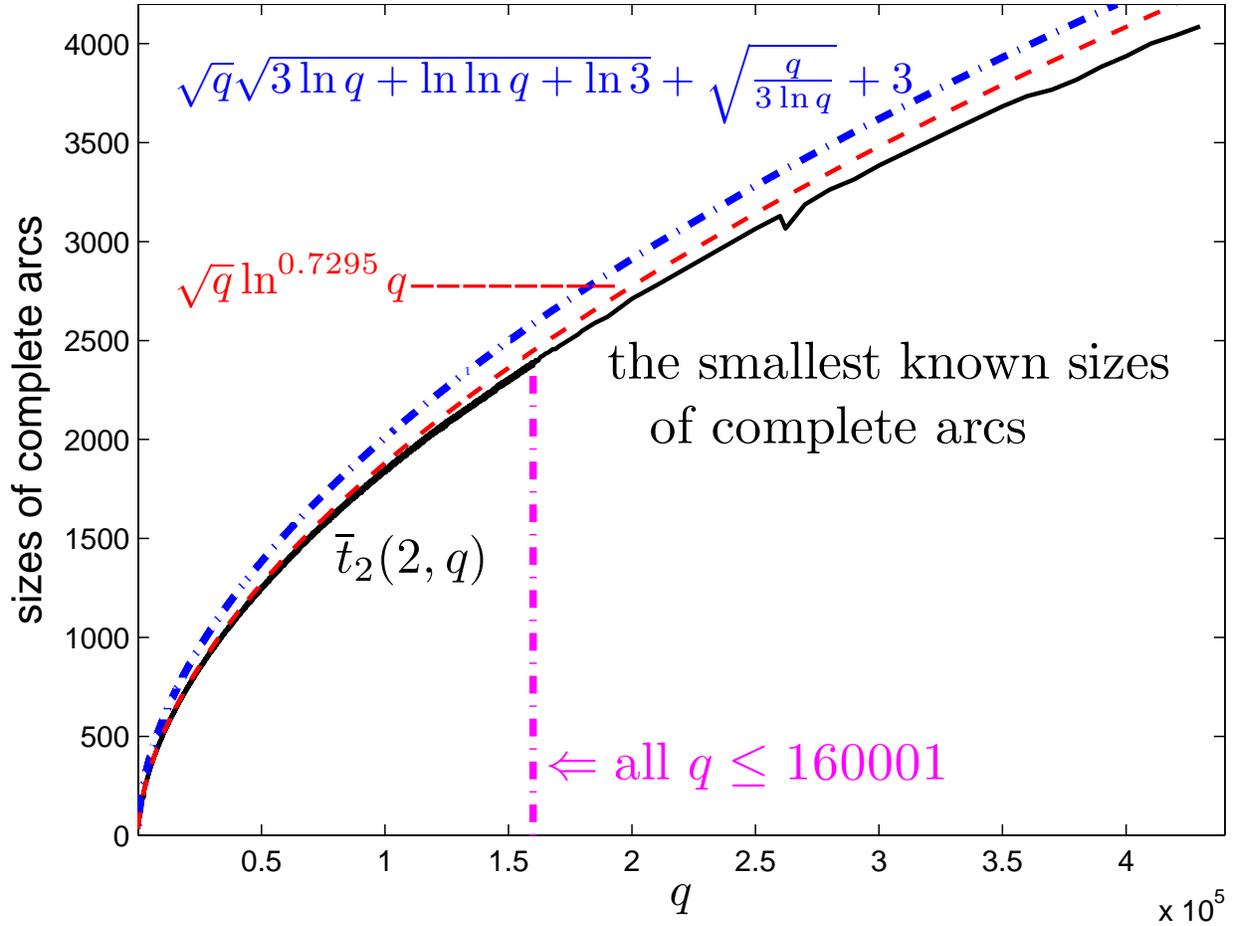}
\caption{\textbf{Upper bounds on $t_{2}(2,q)$ vs the smallest known sizes $\overline{t}_{2}(2,q)$
 of complete arcs.}
Conjectural upper bound from \cite{BDFKMP-PIT2014}
$\sqrt{q}\sqrt{3\ln q+\ln \ln q+\ln
3}+\sqrt{\frac{q}{3\ln q}}+3$ (\emph{top dashed-dotted blue curve}); upper bound $\sqrt{q}\ln^{0.7295}q$
(\emph{the 2-nd dashed red curve});
the smallest known sizes $\overline{t}_{2}(2,q)$ of complete arcs in $\mathrm{PG}(2,q)$, $q\in Q$
(\emph{bottom solid black curve}). \emph{Vertical dashed-dotted magenta line} marks region $q\le 160001$
of the complete search for all $q$'s prime power.}
\label{fig1_t_Gr&bounds}
\end{figure}

 Using sizes $\overline{t}_{2}(2,q)$
for $q\in Q$ collected in Tables 1 -- 6, we obtained the bottom
solid black curve in Figure \ref{fig1_t_Gr&bounds} where the
vertical dashed-dotted magenta line marks the region
$q\le160001$ of complete search for all $q$'s prime power.
Also, in Figure \ref{fig1_t_Gr&bounds} the following upper
bounds on $t_{2}(2,q)$ are shown: the conjectural upper bound
from \cite{BDFKMP-PIT2014} $\sqrt{q}\sqrt{3\ln q+\ln \ln q+\ln
3}+\sqrt{\frac{q}{3\ln q}}+3$ (the top dashed-dotted blue
curve), see \eqref{eq1_Prob bounds}; the bound
$\sqrt{q}\ln^{0.7295}q$ (the 2-nd dashed red curve), see
\eqref{eq1_Bounds2G}.  Note that, on the scale of
$\overline{t}_{2}(2,q)$, the curves of the bounds and of
$\overline{t}_{2}(2,q)$ often almost coalesce with each other,
see Figures \ref{fig4_t_Gr&1.006} and~\ref{fig5_t_Gr&varphi}.
At the same time, non-standard types of upper bounds defined by
\eqref{eq1 c(q)} -- \eqref{eq1_h(q)} allow us to show bounds
more expressive.

\begin{figure}[htbp]
\includegraphics[width=\textwidth]{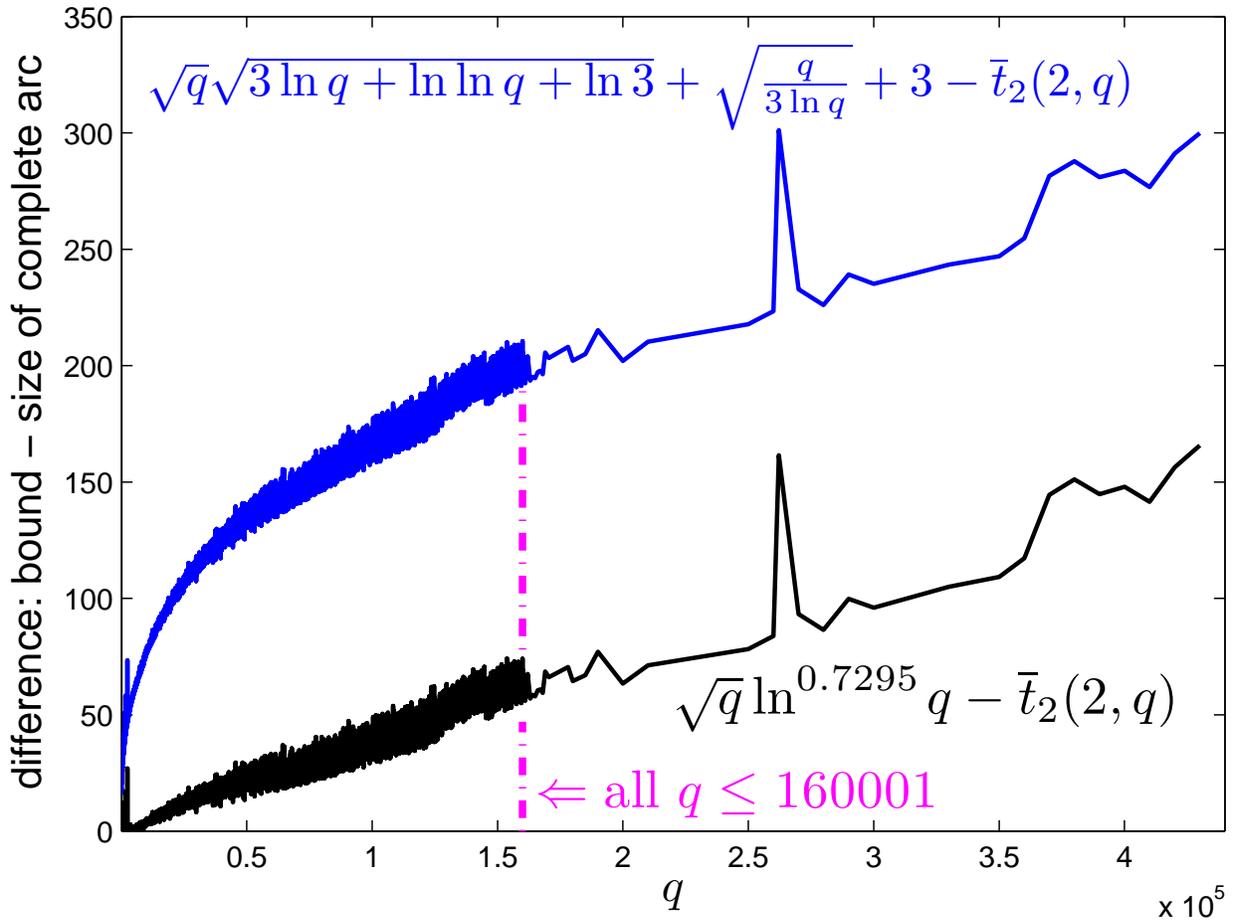}
\caption{\textbf{Differences between upper bounds on $t_{2}(2,q)$
 and the smallest known sizes $\overline{t}_{2}(2,q)$
of complete arcs in $\mathrm{PG}(2,q)$, $q\in Q$ (I):}  differences for conjectural bound
$\sqrt{q}\sqrt{3\ln q+\ln \ln q+\ln
3}+\sqrt{\frac{q}{3\ln q}}+3$ (\emph{top blue curve}); differences for bound $\sqrt{q}\ln^{0.7295} q$
(\emph{bottom black curve}). \emph{Vertical dashed-dotted magenta line} marks region $q\le 160001$
of the complete search for all $q$'s prime power.}
\label{fig2_Delta_conject_0.7295}
\end{figure}

Figure \ref{fig2_Delta_conject_0.7295} presents the differences
between upper bounds on $t_{2}(2,q)$ and the smallest known
sizes $\overline{t}_{2}(2,q)$ of complete arcs in
$\mathrm{PG}(2,q)$, $q\in Q$: the differences
$$\sqrt{q}\sqrt{3\ln q+\ln \ln q+\ln
3}+\sqrt{\frac{q}{3\ln q}}+3 -\overline{t}_{2}(2,q)$$ for the
conjectural upper bound \eqref{eq1_Prob bounds} (the top blue
curve) and the differences
$$\sqrt{q}\ln^{0.7295}q-\overline{t}_{2}(2,q)$$
for the bound \eqref{eq1_Bounds2G} (the bottom black curve).
The vertical dashed-dotted magenta line marks the region
$q\le160001$ of complete search for all $q$'s prime power.

\begin{figure}[htbp]
\includegraphics[width=\textwidth]{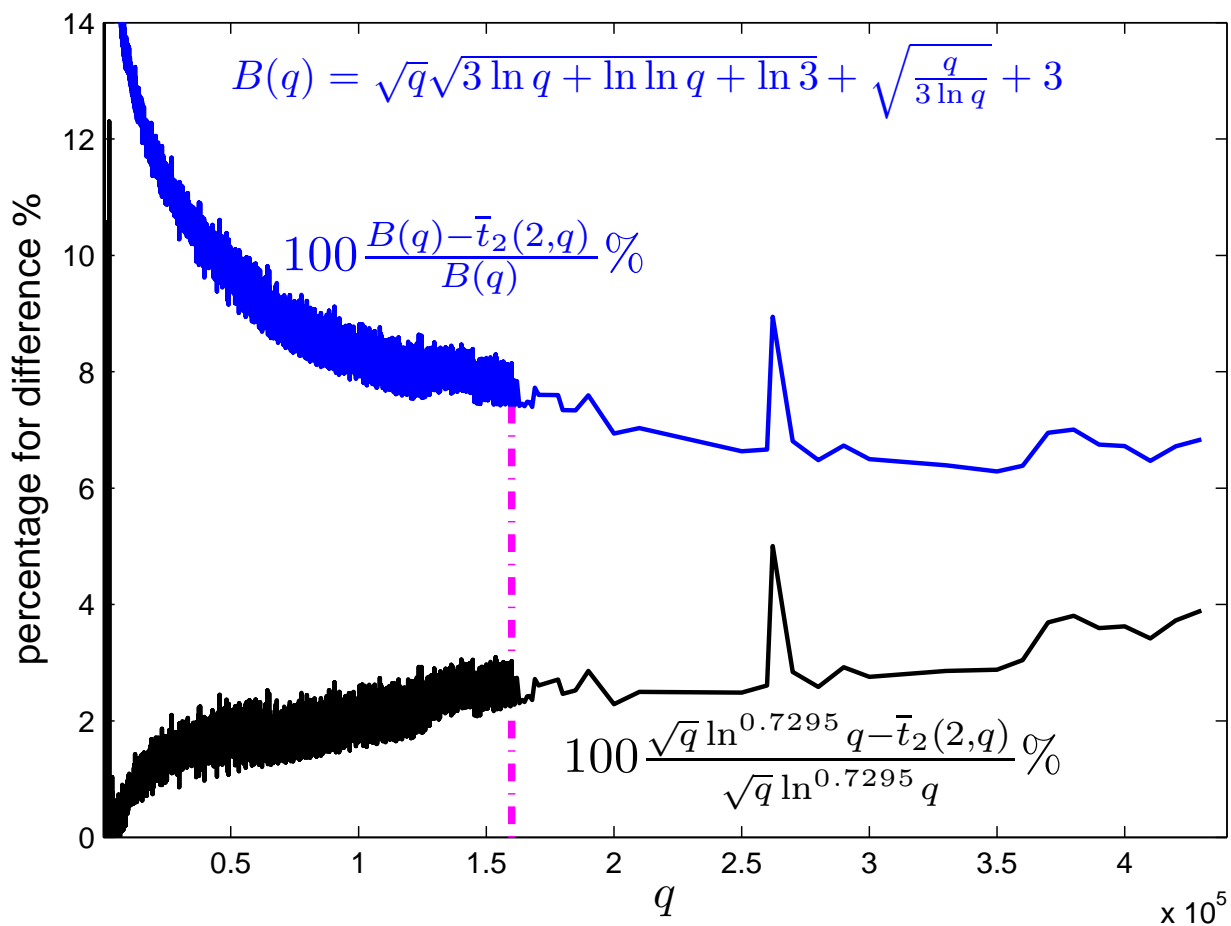}
\caption{\textbf{Percentages for differences between upper bounds on $t_{2}(2,q)$
and the smallest known sizes $\overline{t}_{2}(2,q)$
of complete arcs in $\mathrm{PG}(2,q)$, $q\in Q$ (I):} percentages for conjectural bound
$\sqrt{q}\sqrt{3\ln q+\ln \ln q+\ln
3}+\sqrt{\frac{q}{3\ln q}}+3$ (\emph{top blue curve});
percentages for bound $\sqrt{q}\ln^{0.7295} q$
(\emph{bottom black curve}). \emph{Vertical dashed-dotted magenta line} marks region $q\le 160001$
of the complete search for all $q$'s prime power.}
\label{fig3_percent_conject_0.7295}
\end{figure}

\begin{figure}[htbp]
\includegraphics[width=\textwidth]{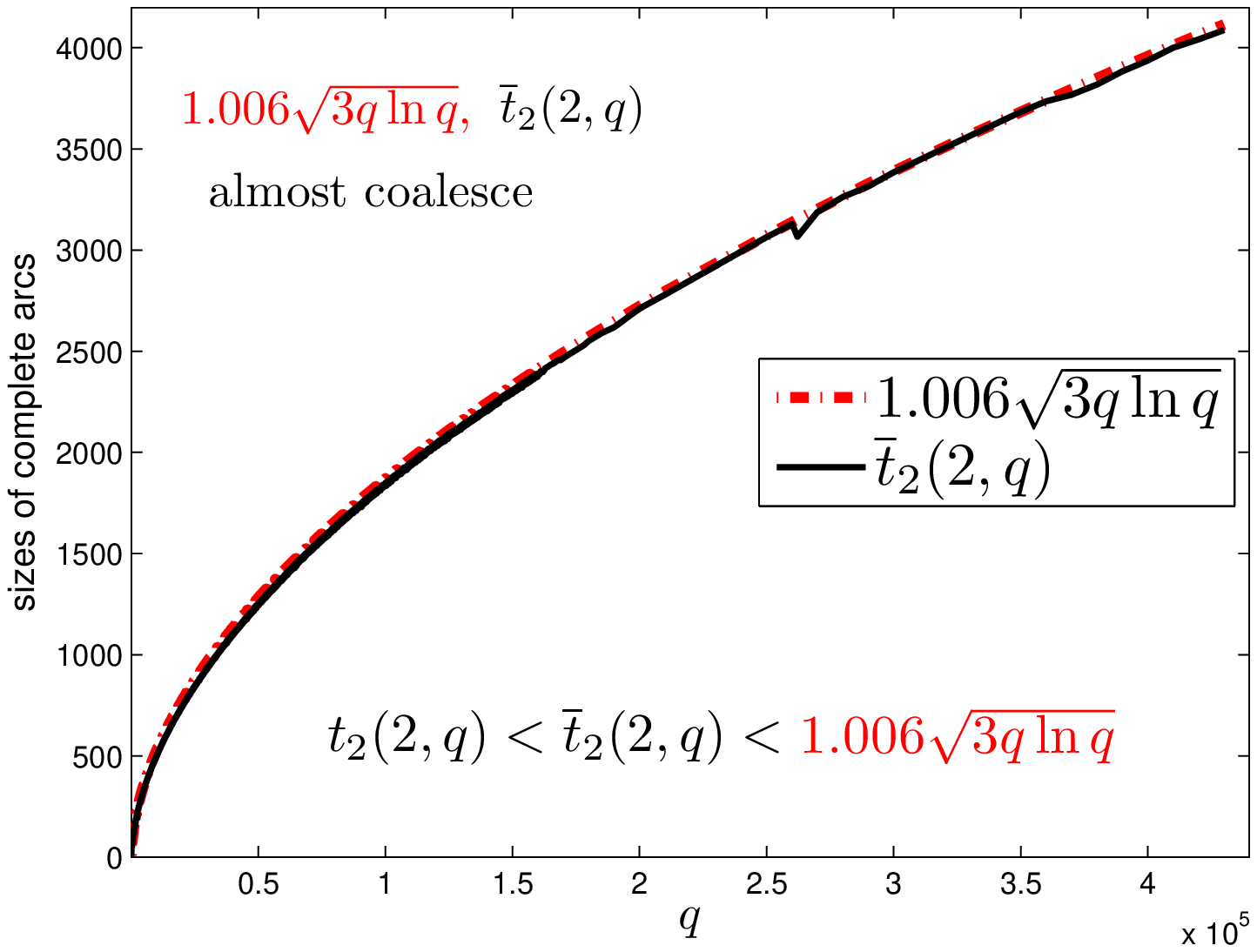}
\caption{\textbf{Upper bound $t_{2}(2,q)<1.006\sqrt{3q\ln q}$ vs the smallest known sizes $\overline{t}_{2}(2,q)$
of complete arcs.}
 Upper bound  $1.006\sqrt{3q\ln q}$
(\emph{top dashed-dotted red curve});
the smallest known sizes $\overline{t}_{2}(2,q)$ of complete arcs in $\mathrm{PG}(2,q)$, $q\in Q$
(\emph{bottom solid black curve}).}
\label{fig4_t_Gr&1.006}
\end{figure}

\begin{figure}[htbp]
\includegraphics[width=\textwidth]{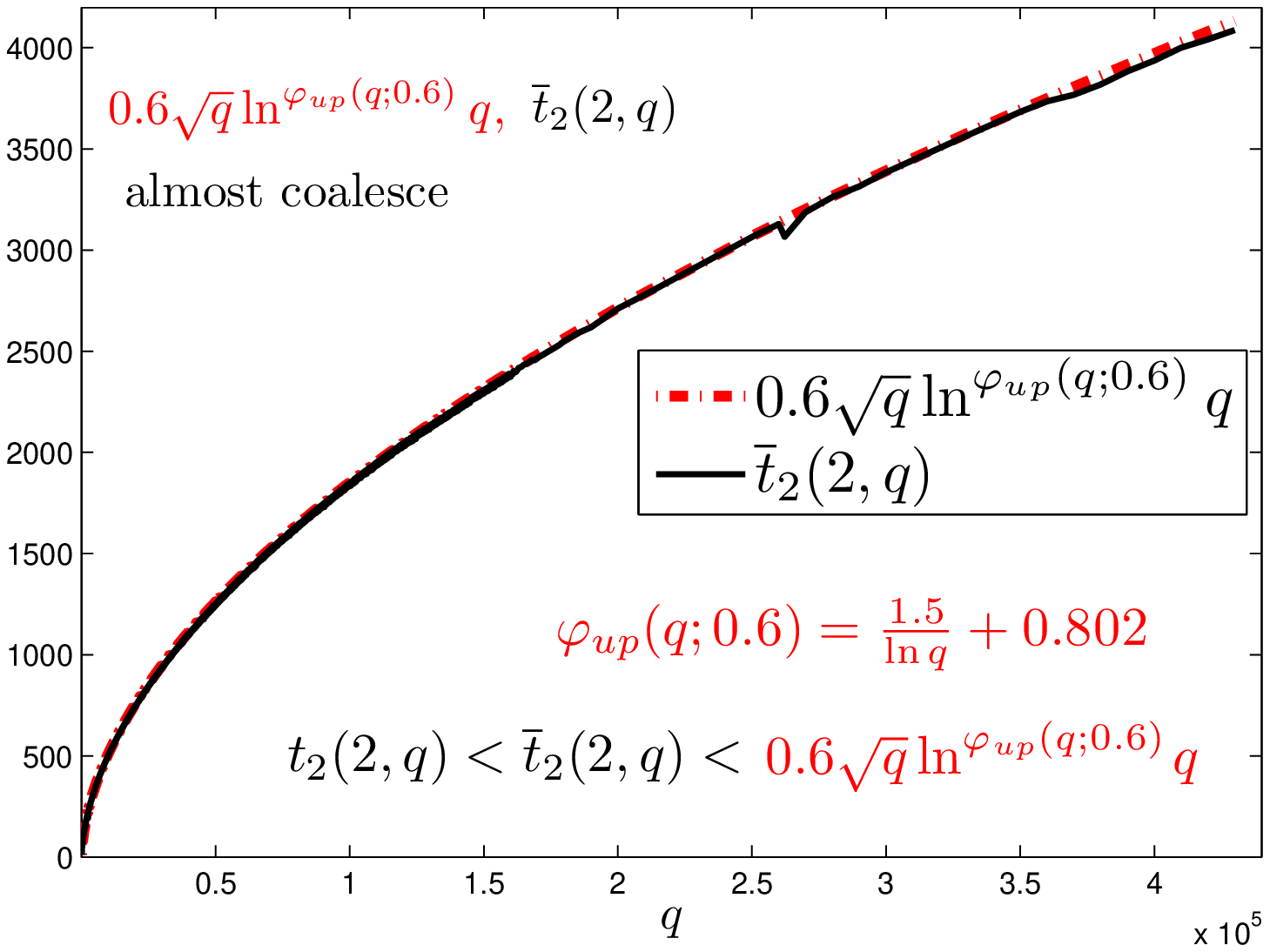}
\caption{\textbf{Upper bound $t_{2}(2,q)<0.6\sqrt{q}\ln ^{\varphi_{up}(q;0.6)}q$
vs the smallest known sizes $\overline{t}_{2}(2,q)$
of complete arcs.}  Upper bound $0.6\sqrt{q}\ln ^{\varphi_{up}(q;0.6)}q$
(\emph{top dashed-dotted red curve});
the smallest known sizes $\overline{t}_{2}(2,q)$ of complete arcs in $\mathrm{PG}(2,q)$, $q\in Q$
(\emph{bottom solid black curve}).}
\label{fig5_t_Gr&varphi}
\end{figure}

Figure \ref{fig3_percent_conject_0.7295} shows the percentages
for differences between upper bounds on $t_{2}(2,q)$ and the
smallest known sizes $\overline{t}_{2}(2,q)$ of complete arcs
in $\mathrm{PG}(2,q)$, $q\in Q$: the percentage
\begin{align*}
100\frac{B(q)-\overline{t}_{2}(2,q)}{B(q)}\%,\qquad
B(q)=\sqrt{q}\sqrt{3\ln q+\ln \ln q+\ln
3}+\sqrt{\frac{q}{3\ln q}}+3,
\end{align*}
for the conjectural upper bound \eqref{eq1_Prob bounds} (the
top blue curve) and the percentage
\begin{align*}
100\frac{\sqrt{q}\ln^{0.7295} q-\overline{t}_{2}(2,q)}{\sqrt{q}\ln^{0.7295} q\vphantom{H^{H^{a}}}}\%
\end{align*}
for the bound  \eqref{eq1_Bounds2G} (the bottom black curve).
The vertical dashed-dotted magenta line marks the region
$q\le160001$ of complete search for all $q$'s prime power.

Figure \ref{fig4_t_Gr&1.006} presents the upper bound
$t_{2}(2,q)<1.006\sqrt{3q\ln q}$ (the top dashed-dotted red
curve), see \eqref{eq1_Bounds1GQ4}, and the smallest known
sizes $\overline{t}_{2}(2,q)$ of complete arcs in
$\mathrm{PG}(2,q)$, $q\in Q$ (the bottom solid black curve).
The curves of the bound $1.006\sqrt{3q\ln q}$ and of
$\overline{t}_{2}(2,q)$ almost coalesce with each other.

We denote
\begin{align}
&c_{up}(q)=\frac{0.27}{\ln q}+0.7\label{eq2_c_up};\\
&\varphi_{up}(q;0.6)=\frac{1.5}{\ln
q}+0.802.\label{eq2_varphi-up}
\end{align}
The bounds \eqref{eq1_Bounds3G} and \eqref{eq1_Bounds4G} are
equivalent, respectively, to
\begin{align}\label{eq2_bound-cup}
&t_{2}(2,q)<\sqrt{q}\ln ^{c_{up}(q)}q;\\
\label{eq2_bound-varphi}
&t_{2}(2,q)<0.6\sqrt{q}\ln ^{\varphi_{up}(q;0.6)}q.
\end{align}

Figure \ref{fig5_t_Gr&varphi} shows the smallest known sizes
$\overline{t}_{2}(2,q)$ of complete arcs in $\mathrm{PG}(2,q)$,
$q\in Q$ (the bottom solid black curve) and the upper bound
$t_{2}(2,q)<0.6\sqrt{q}\ln ^{\varphi_{up}(q;0.6)}q$ (the top
dashed-dotted red curve), see \eqref{eq1_Bounds4G} and
\eqref{eq2_bound-varphi}. The curves of the bound
$0.6\sqrt{q}\ln ^{\varphi_{up}(q;0.6)}q$ and of
$\overline{t}_{2}(2,q)$  almost coalesce with each other.

\begin{figure}[htbp]
\includegraphics[width=\textwidth]{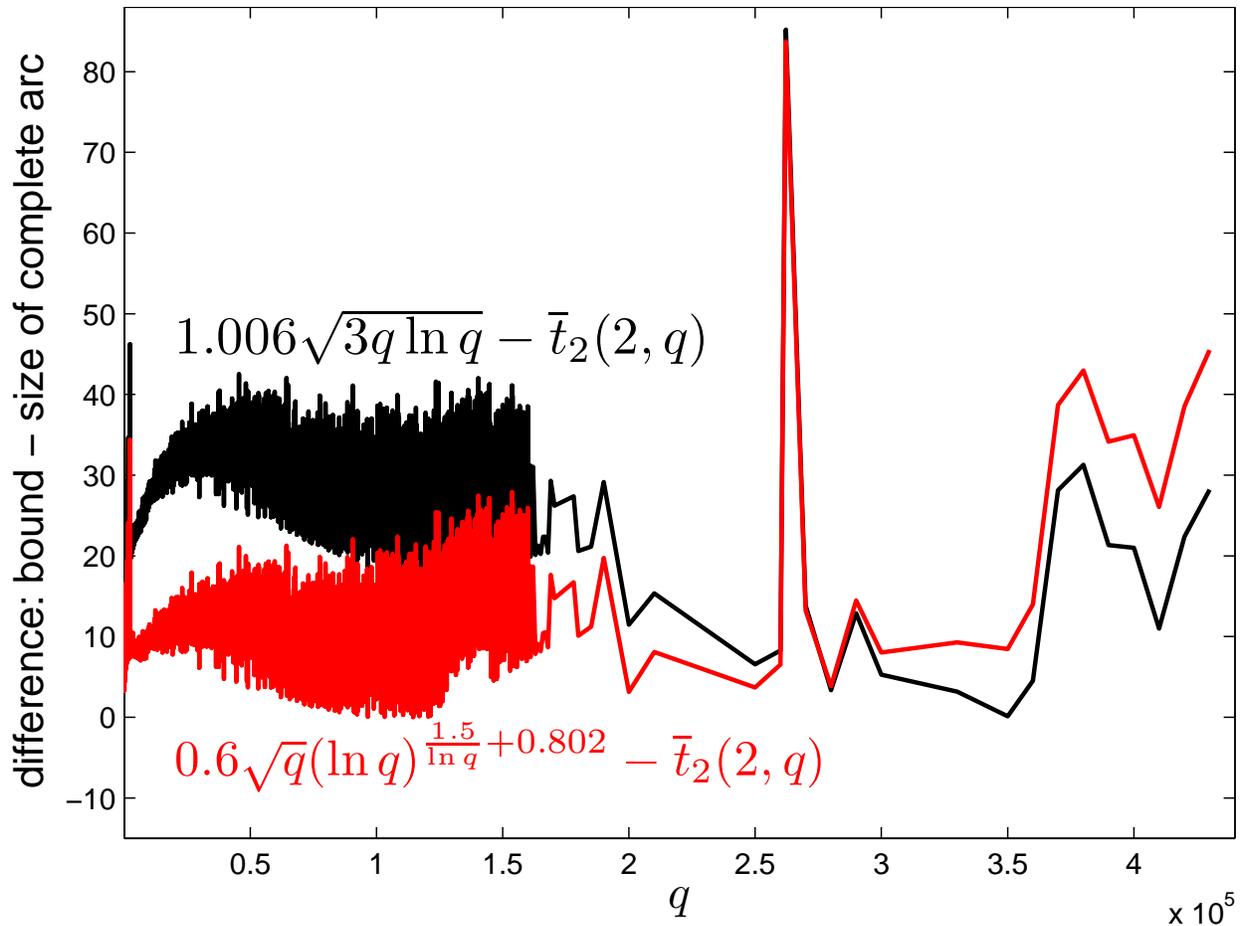}
\caption{\textbf{Differences between upper bounds on $t_{2}(2,q)$
 and the smallest known sizes $\overline{t}_{2}(2,q)$
of complete arcs in $\mathrm{PG}(2,q)$, $q\in Q$ (II):}  differences for bound
$1.006\sqrt{3q\ln q}$ (\emph{top black curve}); differences for bound
$0.6\sqrt{q}\ln^{\varphi_{up}(q;0.6)} q,~\varphi_{up}(q;0.6)=\frac{1.5}{\ln q}+0.802$
(\emph{bottom red curve}).}
\label{fig6_Delta_1.006_0.6}
\end{figure}

Figure \ref{fig6_Delta_1.006_0.6} presents the differences
between upper bounds on $t_{2}(2,q)$ and the smallest known
sizes $\overline{t}_{2}(2,q)$ of complete arcs in
$\mathrm{PG}(2,q)$, $q\in Q$: the differences
$$1.006\sqrt{3q\ln q}-\overline{t}_{2}(2,q)$$
for the bound \eqref{eq1_Bounds1GQ4} (the top black curve) and
the differences
$$0.6\sqrt{q}\ln^{\varphi_{up}(q;0.6)} q-\overline{t}_{2}(2,q),\qquad \varphi_{up}(q;0.6)=\frac{1.5}{\ln q}+0.802,$$
for the bound \eqref{eq1_Bounds4G}, \eqref{eq2_bound-varphi}
(the bottom red curve).

\begin{figure}[htbp]
\includegraphics[width=\textwidth]{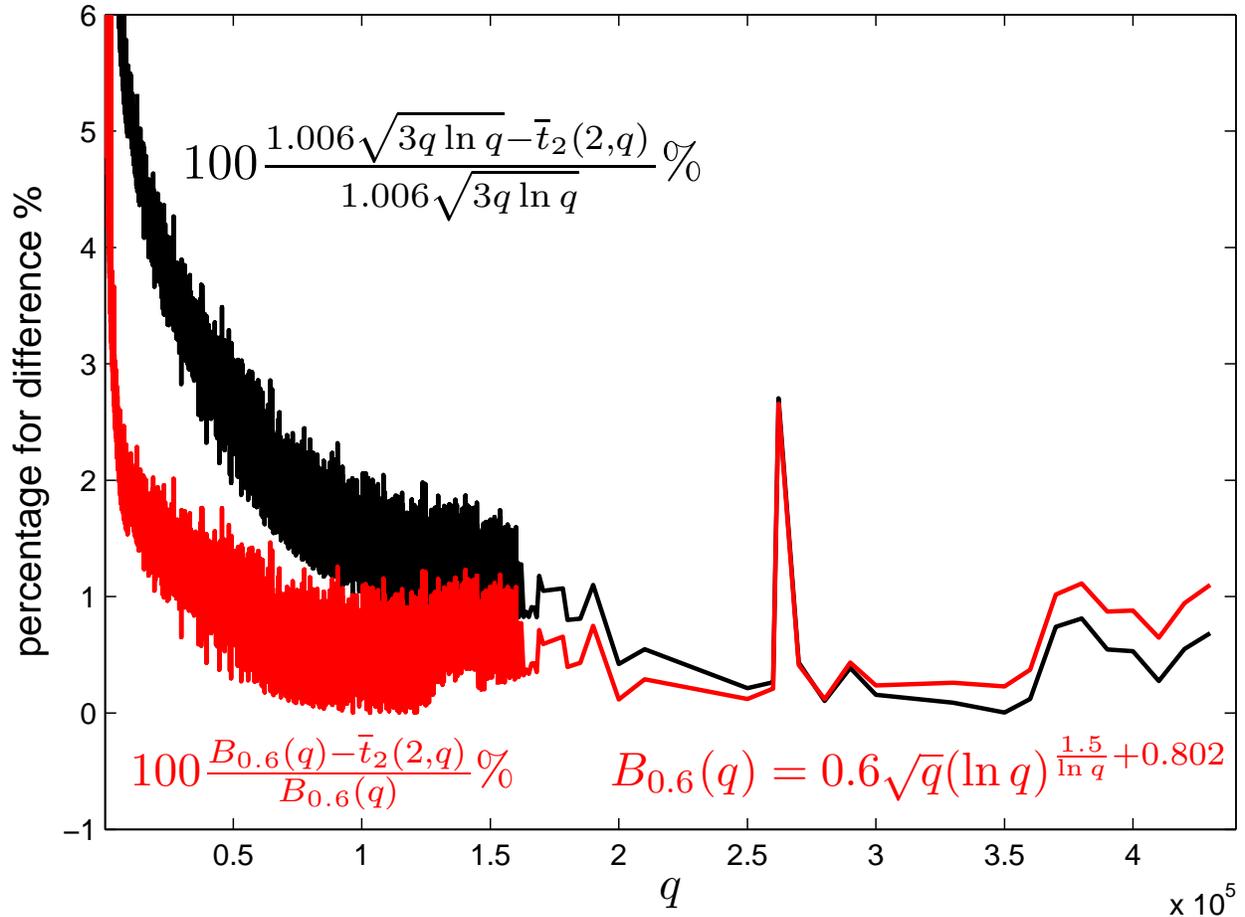}
\caption{\textbf{Percentages for differences between upper bounds on $t_{2}(2,q)$
 and the smallest known sizes $\overline{t}_{2}(2,q)$
of complete arcs in $\mathrm{PG}(2,q)$, $q\in Q$ (II):} percentage for bound
$1.006\sqrt{3q\ln q}$ (\emph{top black curve}); percentage for bound
$0.6\sqrt{q}\ln^{\varphi_{up}(q;0.6)} q,~\varphi_{up}(q;0.6)=\frac{1.5}{\ln q}+0.802$
(\emph{bottom red curve}).}
\label{fig7_percent_1.006_0.6}
\end{figure}

Figure \ref{fig7_percent_1.006_0.6} presents the percentages
for differences between upper bounds on $t_{2}(2,q)$ and the
smallest known sizes $\overline{t}_{2}(2,q)$ of complete arcs
in $\mathrm{PG}(2,q)$, $q\in Q$: the percentage
\begin{align*}
100\frac{1.006\sqrt{3q\ln q}-\overline{t}_{2}(2,q)}{1.006\sqrt{3q\ln q}}\%
\end{align*}
for the upper bound \eqref{eq1_Bounds1GQ4} (the top black
curve) and the percentage
\begin{align*}
100\frac{B_{0.6}(q)-\overline{t}_{2}(2,q)}{B_{0.6}(q)}\%,
~B_{0.6}(q)=0.6\sqrt{q}\ln^{\varphi_{up}(q;0.6)} q,~\varphi_{up}(q;0.6)=\frac{1.5}{\ln q}+0.802,
\end{align*}
for the bound  \eqref{eq1_Bounds4G}, \eqref{eq2_bound-varphi}
(the bottom red curve).

\begin{figure}[htbp]
\includegraphics[width=\textwidth]{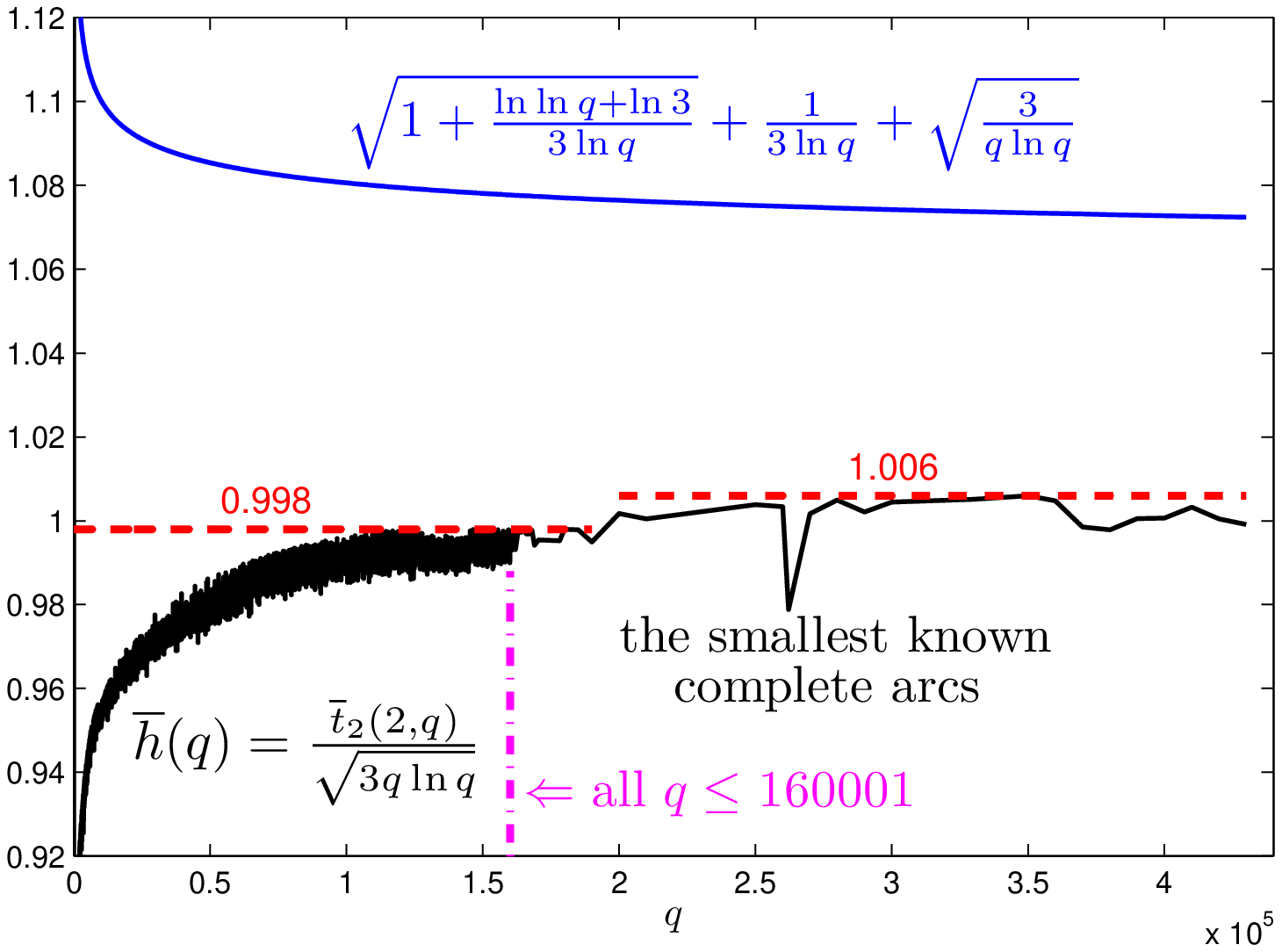}
\caption{Conjectural upper bound $\sqrt{1+(\ln\ln q+\ln3)/(3\ln q)}+1/(3\ln q)+\sqrt{3/(q\ln q)}$
 (top solid blue curve); values of
$\overline{h}(q)=\overline{t}_{2}(2,q)/\sqrt{3q\ln q}$\,
for the smallest known complete arcs in $\mathrm{PG}(2,q)$, $q\in Q$ (\emph{bottom solid black curve});
upper bounds ``0.998'' and ``1.006'' (\emph{horizontal dashed red lines}).
\emph{Vertical dashed-dotted magenta line} marks region $q\le 160001$
of the complete search for all $q$'s prime power.}
\label{fig8_h_Gr&bounds}
\end{figure}

Figure \ref{fig8_h_Gr&bounds} shows the conjectural upper bound
\begin{align}\label{eq2_bound-conject_divided}
\sqrt{1+\frac{\ln\ln q+\ln3}{3\ln q}}+\frac{1}{3\ln q}+\sqrt{\frac{3}{q\ln q}}
\end{align}
obtained by division of the bound \eqref{eq1_Prob bounds} by
$\sqrt{3q\ln q}$ (the top solid blue curve). Also, values of
$\overline{h}(q)=\overline{t}_{2}(2,q)/\sqrt{3q\ln q}$ for the
smallest known complete arcs in $\mathrm{PG}(2,q)$, $q\in Q$,
are presented (the bottom solid black curve). The bounds
\eqref{eq1_Bounds1G} and \eqref{eq1_Bounds1GQ4} are shown by
the dashed red lines $y=0.998$ and $y=1.006$. The vertical
dashed-dotted magenta line marks the region $q\le 160001$ of
complete search for all $q$'s prime power.

\begin{figure}[htbp]
\includegraphics[width=\textwidth]{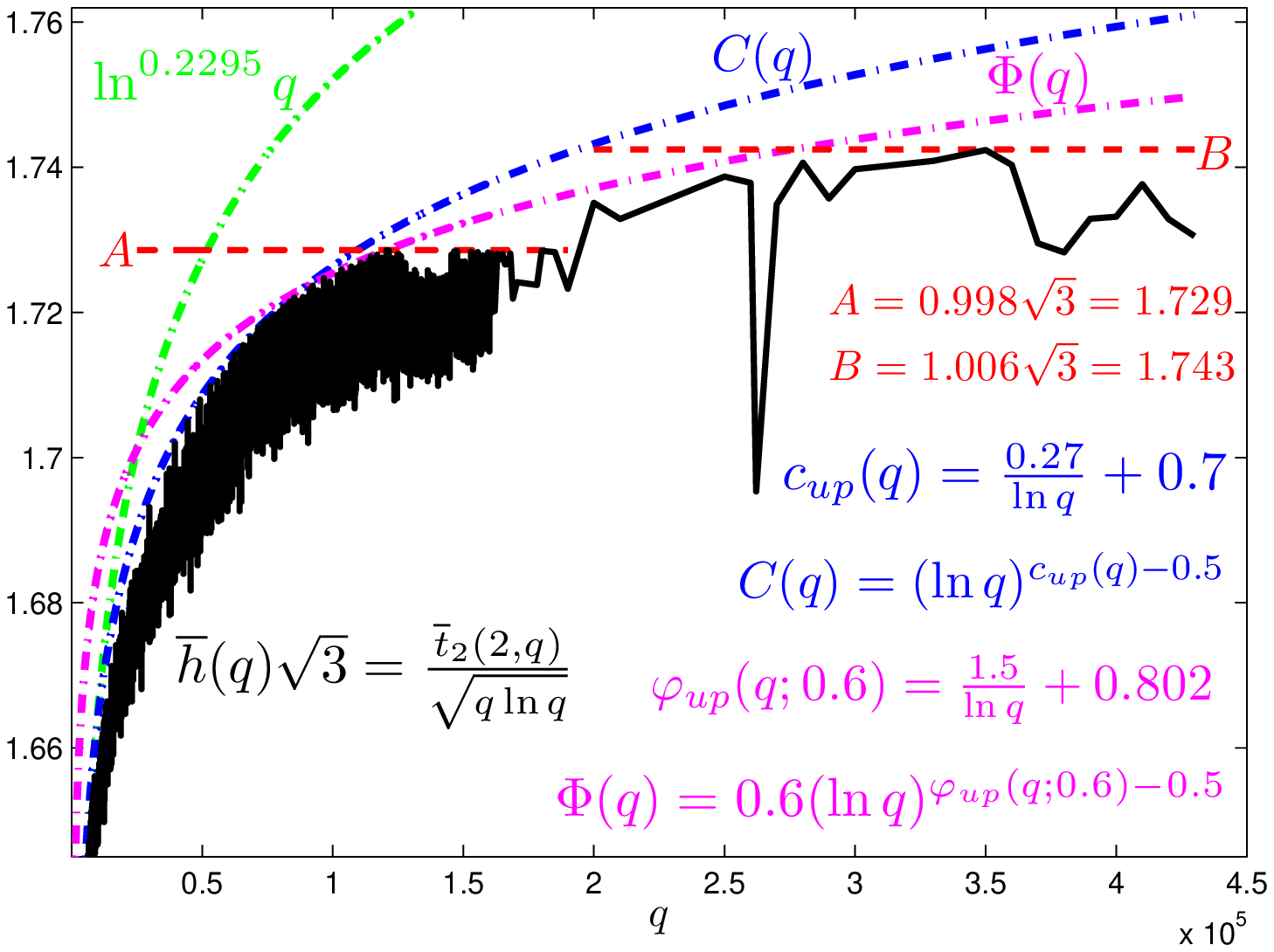}
\caption{Values of
$\overline{h}(q)\sqrt{3}=\overline{t}_{2}(2,q)/\sqrt{q\ln q}$ for
the smallest known complete arcs in $\mathrm{PG}(2,q)$, $q\in Q$ (\emph{bottom solid black curve}) and
upper bounds on $t_{2}(2,q)/\sqrt{q\ln q}$:\quad
 $\ln^{0.2295}q$ (\emph{top dashed-dotted green curve}),
$C(q)=(\ln q)^{c_{up}(q)-0.5} $ with $c_{up}(q)
=\frac{0.27}{\ln q}+0.7$ (the 2-nd dashed-dotted blue curve),
$\Phi(q)=0.6(\ln q)^{\varphi_{up}(q;0.6)-0.5}$ with
$\varphi_{up}(q;0.6)=\frac{1.5}{\ln q}+0.802$ (\emph{the 3-rd dashed-dotted magenta curve}),
$0.998\sqrt{3}=1.729$ (\emph{dashed red line} $A$), $1.006\sqrt{3}=1.743$ (\emph{dashed red line} $B$)}
\label{fig9_hsqrt(3)_Gr&bounds}
\end{figure}

Figure \ref{fig9_hsqrt(3)_Gr&bounds} shows values of
$\overline{h}(q)\sqrt{3}=\overline{t}_{2}(2,q)/\sqrt{q\ln q}$
for the smallest known complete arcs in $\mathrm{PG}(2,q)$,
$q\in Q$ (the bottom solid black curve) and the following upper
bounds on $t_{2}(2,q)/\sqrt{q\ln q}$:\quad $\ln^{0.2295}q$ (the
top dashed-dotted green curve), $C(q)=(\ln q)^{c_{up}(q)-0.5}$
with $c_{up}(q)=\frac{0.27}{\ln q}+0.7$ (the 2-nd dashed-dotted
blue curve), $\Phi(q)=0.6(\ln q)^{\varphi_{up}(q;0.6)-0.5}$
with $\varphi_{up}(q;0.6)=\frac{1.5}{\ln q}+0.802$ (the 3-rd
dashed-dotted magenta curve), $0.998\sqrt{3}=1.729$ (the dashed
red line $A$), $1.006\sqrt{3}=1.743$ (the dashed red line $B$).
Correlations of these bounds are given by the following
theorem.

\begin{theorem}\label{th2_minimum}
Let $t_{2}(2,q)$ be the smallest size of a complete arc in the
projective plane $\mathrm{PG}(2,q)$. Let $Q$ be the set of
values of $q$ given by \eqref{eq1_Q}. For $q\in Q$, in
$\mathrm{PG}(2,q)$ the following upper bounds hold:
    \begin{align}
&t_{2}(2,q)<
\min\{\sqrt{q}\ln^{0.7295}q,~\sqrt{q}\ln^{c_{up}(q)} q,~0.6\sqrt{q}
\ln^{\varphi_{up}(q;0.6)} q\}=\bigskip\bigskip
 \label{eq2_Bounds-min}\displaybreak[3]\\
 &\left\{\begin{array}{lcl}
 \phantom{0.6}\sqrt{q}\ln^{0.7295}q &\text{if}&\phantom{11}109\le q\le 9437\\
 \phantom{0.6}\sqrt{q}\ln^{c_{up}(q)} q &\text{if}&9437< q\le 88873\\
0.6\sqrt{q}\ln^{\varphi_{up}(q;0.6)} q&\text{if}&88873< q
\end{array}
\right.,\notag\\
&\mbox{where }c_{up}(q)=\frac{0.27}{\ln q}+0.7,~~
\varphi_{up}(q;0.6)=\frac{1.5}{\ln
q}+0.802.\notag
\end{align}
\end{theorem}

\begin{figure}[htbp]
\includegraphics[width=\textwidth]{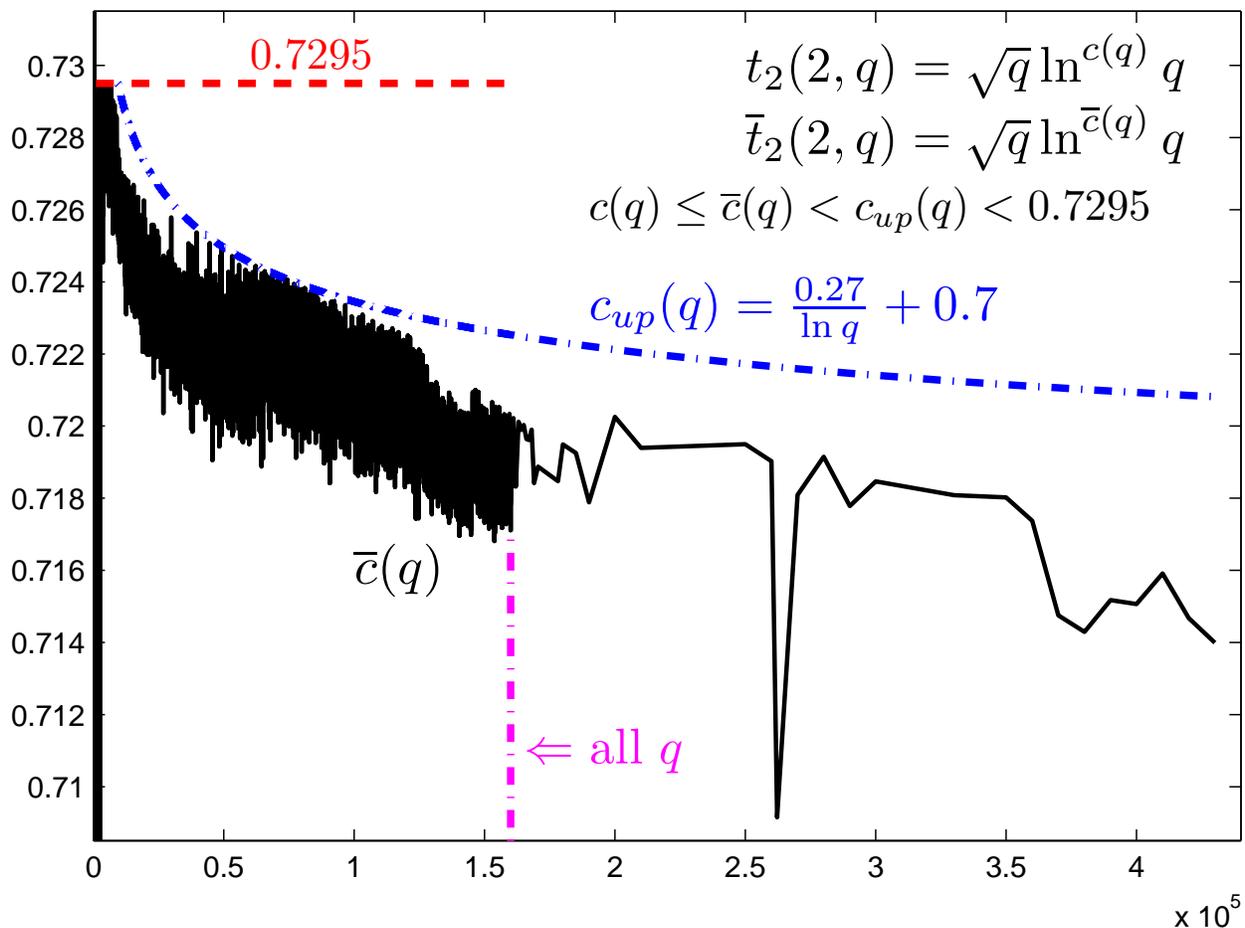}
\caption{\textbf{Values and bounds connected with decreasing function $c(q)$:}
 upper bound``0.7295'' (\emph{dashed red line});
  upper bound
 $c_{up}(q)=\frac{0.27}{\ln q}+0.7$ (\emph{top dashed-dotted blue curve});
values of $\overline{c}(q)$ for the smallest known complete
arcs in $\mathrm{PG}(2,q)$, $q\in Q$ (\emph{bottom solid black curve}).
\emph{Vertical dashed-dotted magenta line} marks region $q\le 160001$
of the complete search for all $q$'s prime power.}
\label{fig10_c_Gr&bounds}
\end{figure}

We denote
\begin{align*}
\overline{t}_{2}(2,q)=\sqrt{q}\ln^{\overline{c}(q)} q.
\end{align*}
By \eqref{eq1_main-ineq}, \eqref{eq1 c(q)},
\eqref{eq2_Bounds-min}, we have
\begin{align*}
c(q)\le\overline{c}(q)<\min\{0.7295,~c_{up}(q)\},~~ q\in Q.
\end{align*}
Using  Figure \ref{fig9_hsqrt(3)_Gr&bounds} and Theorem
\ref{th2_minimum}, we obtain Figure \ref{fig10_c_Gr&bounds}
where values and bounds connected with the decreasing function
$c(q)$ of \eqref{eq1 c(q)} are presented. Values of
$\overline{c}(q)$ for the smallest known complete arcs in
$\mathrm{PG}(2,q)$, $q\in Q$, are shown by the bottom solid
black curve. The upper bound \eqref{eq1_Bounds2G} is given by
the dashed red line $y=0.7295$. Finally, the bound
\eqref{eq1_Bounds3G} is shown by the top dashed-dotted blue
curve $c_{up}(q)=\frac{0.27}{\ln q}+0.7$. The vertical
dashed-dotted magenta line marks the region $q\le160001$ of
complete search for all $q$'s prime power.

\begin{figure}[htbp]
\includegraphics[width=\textwidth]{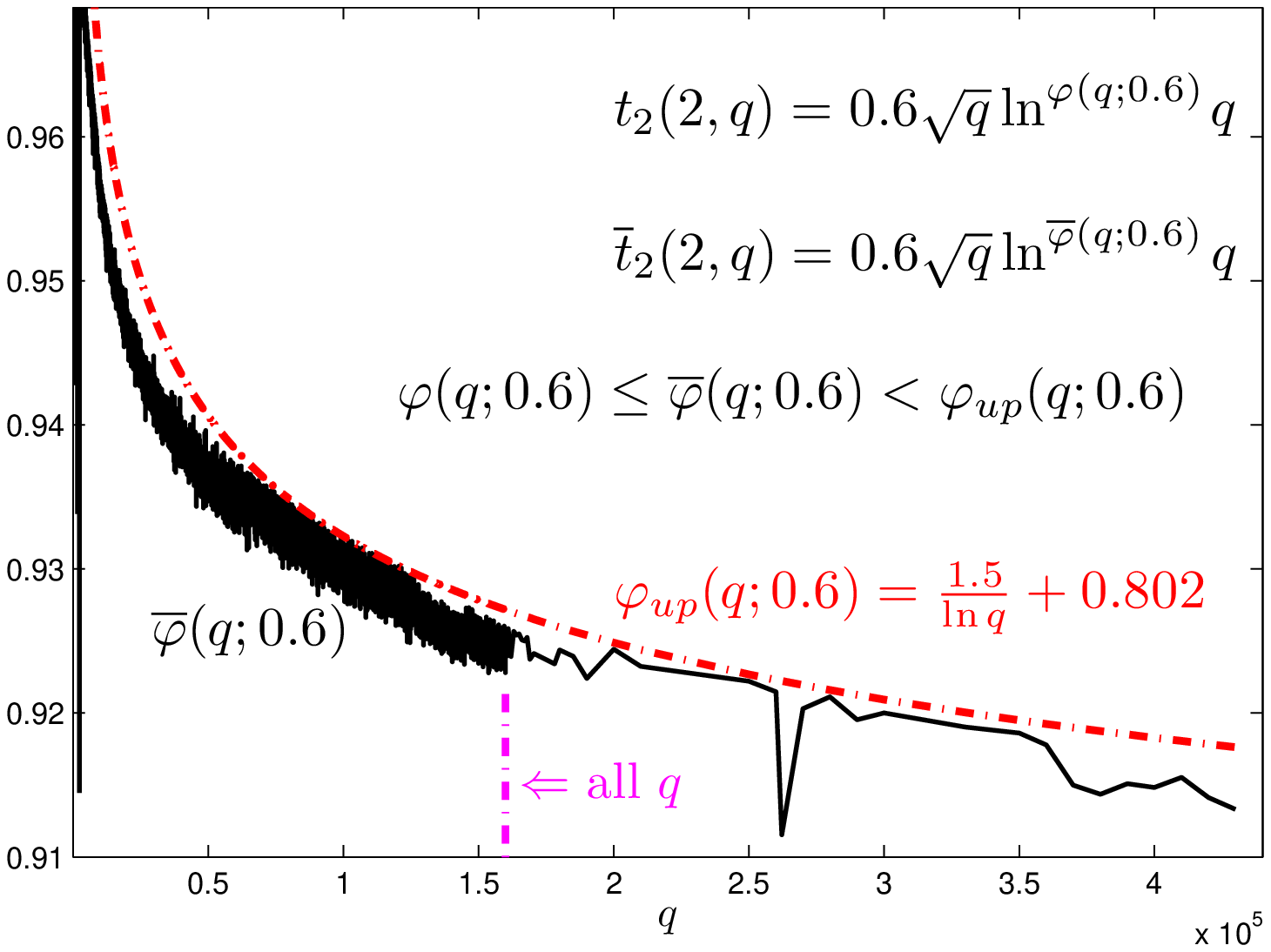}
\caption{\textbf{Values and bounds connected with decreasing function $\varphi(q;0.6)$:}
   upper bound
 $\varphi_{up}(q;0.6)=\frac{1.5}{\ln q}+0.802$ (\emph{top dashed-dotted red curve});
values of $\overline{\varphi}(q;0.6)$ for the smallest known complete
arcs in $\mathrm{PG}(2,q)$, $q\in Q$ (\emph{bottom solid black curve}).
\emph{Vertical dashed-dotted magenta line} marks region $q\le 160001$
of the complete search for all $q$'s prime power.}
\label{fig11_varphi_Gr&bounds}
\end{figure}

We denote
\begin{align*}
\overline{t}_{2}(2,q)=0.6\sqrt{q}\ln^{\overline{\varphi}(q;0.6)} q.
\end{align*}
Figure \ref{fig11_varphi_Gr&bounds} shows the upper bound
 $\varphi_{up}(q;0.6)=\frac{1.5}{\ln q}+0.802$ (the top dashed-dotted red
 curve) and
values of $\overline{\varphi}(q;0.6)$ for the smallest known
complete arcs in $\mathrm{PG}(2,q)$, $q\in Q$ (the bottom solid
black curve). The vertical dashed-dotted magenta line marks the
region $q\le160001$ of complete search for all $q$'s prime
power.

\begin{remark}
\label{rem2_tooth} In Figures
\ref{fig1_t_Gr&bounds}--\ref{fig11_varphi_Gr&bounds}, the
``tooth'' connected with $\overline{t}_{2}(2,2^{18})=3066$ is
provided by the construction of \cite{DGMP-Innov}. It reminds
to the reader that, in principle, theoretical constructions are
able to give a better result than greedy algorithms. On the
other side, as a rule, the orders of arc size from a
theoretical construction and a greedy algorithm are similar. We
have $2^{18}=262144$; the nearest known greedy algorithm's
result is $\overline{t}_{2}(2,260003)=3129$, see Table~3. The
difference between 3129 and 3066 is $\thickapprox 2\%$.
\end{remark}

\section{List of tables}

\hspace{0.5cm}\textbf{Table 1.} Sizes
$\overline{t}_{2}=\overline{t} _{2}(2,q)$ of complete arcs in
$\mathrm{PG}(2,q)$ smaller than in \cite{BDFMP-JG2013}.
\textbf{p. 30}

\textbf{Table 2.} The smallest known sizes
 $\overline{t}_{2}=\overline{t}_{2}(2,q)$
 of complete arcs in planes
 $\mathrm{PG}(2,q),$ $49729\leq q\leq 160801$, $q$
 non-prime. \textbf{p. 30}

 \textbf{Table 3.} The smallest known sizes
 $\overline{t}_{2}=\overline{t}_{2}(2,q)$
 of complete arcs in planes
 $\mathrm{PG}(2,q),$ $q=124477,128683,160001$, and $160801\leq q\leq 430007$, $q$
 sporadic. \textbf{p. 30}

\textbf{Table 4.} The smallest known sizes
 $\overline{t}_{2}=\overline{t}_{2}(2,q)$
 of complete arcs in
 $\mathrm{PG}(2,q),$\\ $2\leq q\leq 10000$, $q$  power prime. \textbf{pp. 31--34}

\textbf{Table 5.} The smallest known sizes
 $\overline{t}_{2}=\overline{t}_{2}(2,q)$
 of complete arcs in planes
 $\mathrm{PG}(2,q),$ $10001\leq q\leq 100000$, $q$  power prime.
 \textbf{pp. 35--64}

 \textbf{Table 6.} The smallest known sizes
 $\overline{t}_{2}=\overline{t}_{2}(2,q)$
 of complete arcs in planes
 $\mathrm{PG}(2,q),$ $100001\leq q\leq 160000$, $q$  power prime.
\textbf{ pp. 65--85}

\newpage
\section*{Appendix. Tables of the smallest known sizes $\overline{t} _{2}(2,q)$
of complete arcs in the plane $\mathrm{PG}(2,q)$}

\textbf{Table 1.} Sizes $\overline{t}_{2}=\overline{t}
_{2}(2,q)$ of complete arcs in planes $\mathrm{PG}(2,q)$
smaller than in \cite{BDFMP-JG2013}\footnote{for $q\le151$ the
sizes are obtained in \cite{Pace-A5}; for $q\ge2383$ the sizes
are obtained in \cite{BDFKMP-ArXiv2013}}
\medskip\\
{\small
\renewcommand{\arraystretch}{1.0}
\begin{tabular}{@{}c@{\,\,\,\,}c|c@{\,\,\,\,}c|c@{\,\,\,\,}c|c@{\,\,\,\,}c|
c@{\,\,\,\,}c|c@{\,\,\,\,}c|c@{\,\,\,\,}c|c@{\,\,\,\,}c@{}}
\hline
$q$ & $\overline{t}_{2}$ & $q$ & $\overline{t}_{2}$ & $q$ & $\overline{t}_{2}$ & $q$ & $
\overline{t}_{2}$ & $q$ & $\overline{t}_{2}$ & $q$ & $\overline{t}_{2}$& $q$
&$\overline{t}_{2}$ & $q$ &$\overline{t}_{2}^{\vphantom{H^{L^{L}}}}$
\\ \hline
67 & 22  &97 & 28  & 109 & 30  &
121 & 30  & 125 & 30  & 137 & 36 &139 & 36  & 151 & 36\\ \hline
2383 & 217 & 2551 & 226 & 2803 & 239 & 2843 & 241 & 2861 & 241 & 2963 & 247& 3125 & 255 & 3167 & 257
\\
 3209 & 259 & 3251 & 261 & 3313 & 263 & 3463 & 271 & 3529 & 274 & 3571 & 274 & 3593 & 277 & 3637 & 279
\\
  3659 & 280 & 3727 & 283 &  3769 & 285 & 3793 & 286 & 3863 & 289 & 3907 & 291 & 3929 & 292 & 3931 & 292
\\
4091 & 299 & 4139 & 301 & 4231 & 305 & 4327 & 309 & 4349 & 310 & 4373 & 311 &  4397 & 312 & 4421 & 313
\\
 4493 & 315 & 4519 & 317 & 4567 & 319 & 4591 & 320& 4639 & 322 & 4787 & 328 & 4789 & 327 & 4813 & 329
\\
  4889 & 331 & 4913 & 333&4937 & 333 & 5167 & 343 & 5297 & 347 & 5503 & 356 & 5531 & 356 & 5639 & 361
\\
5827 & 368 & 5881 & 370 & 6211 & 382 & 6491 & 391 &  &  &  &  \\ \hline
\end{tabular}
}

\begin{center}
\textbf{Table 2.} The smallest known sizes
 $\overline{t}_{2}=\overline{t}_{2}(2,q)$
 of complete arcs in planes
 $\mathrm{PG}(2,q),$ $49729\leq q\leq 160801$, $q$
 non-prime.\footnote{The sizes are obtained in \cite{BDFKMP-PIT2014} and \cite{BDFKMP-JGtoappear}.
 All the sizes are new in comparison with \cite{BDFKMP-ArXiv}.}
\medskip\\{\small
 \renewcommand{\arraystretch}{1.00}
 \begin{tabular}
 {rc|rc|rc|rc|rc}
  \hline
$q\quad$ & $\overline{t}_{2}$ &
$q\quad$ & $\overline{t}_{2}$ &
$q\quad$ & $\overline{t}_{2}$ &
$q\quad$ & $\overline{t}_{2}$ &
$q\quad$ & $\overline{t}_{2}^{ \vphantom{H^{L}}}$ \\
 \hline
 49729 & 1246  &  50653 & 1262  &  51529 & 1272  &  52441 & 1284  &  54289 & 1311  \\
 57121 & 1349  &  58081 & 1363  &  59049 & 1370  &  63001 & 1417  &  65536 & 1453  \\
 66049 & 1462  &  68921 & 1497  &  69169 & 1497  &  72361 & 1534  &  73441 & 1552  \\
 76729 & 1589  &  78125 & 1606  &  78961 & 1614  &  79507 & 1623  &  80089 & 1624  \\
 83521 & 1665  &  85849 & 1691  &  94249 & 1785  &  96721 & 1810  &  97969 & 1820  \\
 100489 & 1847  &  103823 & 1878  &  109561 & 1942  &  113569 & 1976  &  117649 & 2011  \\
 120409 & 2036  &  121801 & 2055  &  124609 & 2077  &  128881 & 2116  &  130321 & 2130  \\
 131072 & 2141  &  134689 & 2163  &  139129 & 2208  &  143641 & 2253  &  146689 & 2275  \\
 148877 & 2297  &  151321 & 2312  &  157609 & 2371  &  160801 & 2388  &  &  \\
 \hline
 \end{tabular}
}
\end{center}

\textbf{Table 3.} The smallest known sizes
 $\overline{t}_{2}=\overline{t}_{2}(2,q)$
 of complete arcs in planes
 $\mathrm{PG}(2,q),$ $q=124477,128683,160001$, and $160801\leq q\leq 430007$, $q$
 sporadic.\footnote{The sizes  that are new or smaller in comparison with \cite{BDFKMP-ArXiv}
 are written either in bold font (if they are obtained in \cite{BDFKMP-JGtoappear}) or in italic-bold
 one (if they are obtained in this work).}
\medskip\\\small{
\renewcommand{\arraystretch}{1.00}
 \begin{tabular}
 {@{}c@{\,\,}c@{\,}|@{\,}c@{\,\,}c@{\,}|c@{\,}c@{\,}|@{\,}c@{\,\,}c@{\,}|@{\,}c@{\,\,}c@{\,}|@{\,}c@{\,\,}c@{}}
  \hline
$q$ & $\overline{t}_{2}$ &
$q$ & $\overline{t}_{2}$ &
$q$ & $\overline{t}_{2}$ &
$q$ & $\overline{t}_{2}$ &
$q$ & $\overline{t}_{2}$ &
$q$ & $\overline{t}_{2}^{ \vphantom{H^{L}}}$ \\
 \hline
 \textbf{124477} & \textbf{2087}  &  \textbf{128683} & \textbf{2117}  &  \textbf{160001} &\textbf{ 2383}
 &  \textbf{160801} & \textbf{2388}  &  \textbf{161009} & \textbf{2401}  &  \textbf{162007} & \textbf{2398}  \\
 \emph{\textbf{163003}} & \emph{\textbf{2417}}  &  \textbf{164011} & \textbf{2424}  &  \textbf{165001} & \textbf{2433}
 &  \emph{\textbf{166013} }& \emph{\textbf{2439}}  &  \emph{\textbf{167009}} & \emph{\textbf{2447}}  &  \emph{\textbf{168013}} & \emph{\textbf{2457}}\\
 \emph{\textbf{169003}} & \emph{\textbf{2456}}  &  \emph{\textbf{170503}} & \emph{\textbf{2471}}  &  \textbf{178169} & \textbf{2530}
 &  \emph{\textbf{180001}} & \emph{\textbf{2551}}  &  \emph{\textbf{185021}} & \emph{\textbf{2589} } &  \emph{\textbf{190027} }& \emph{\textbf{2619}}  \\
 200003 & 2711  &  210011 & 2780  &  250007 & 3065  &  260003 & 3129  &  \textbf{262144} & \textbf{3066}  &  270001 & 3188  \\
 280001 & 3262  &  \emph{\textbf{290011}} & \emph{\textbf{3315}}  &  300007 & 3384  &  \textbf{330017}
 & \textbf{3565}  &  \emph{\textbf{350003}} & \emph{\textbf{3683}}  &  \emph{\textbf{360007}} & \emph{\textbf{3735}}  \\
 370003 & 3767  &  380041 & 3819  &  390001 & 3883  &  400009 & 3937  &  410009 & 4000  &  \textbf{420001} & \textbf{4041}  \\
 \textbf{430007} & \textbf{4087}  &  &  &  &  &  &  &  &  &  &  \\
 \hline
 \end{tabular}
 }

 \textbf{Table 4.} The smallest known sizes
 $\overline{t}_{2}=\overline{t}_{2}(2,q)$
 of complete arcs in planes
 $\mathrm{PG}(2,q),$\\ $2\leq q\leq 10000$, $q$ power prime \medskip\\\small{
  \renewcommand{\arraystretch}{0.96}
 \begin{tabular}
 {@{}r@{\,\,\,\,\,}r@{\,\,}|@{\,\,}r@{\,\,\,\,\,}r@{\,\,}|@{\,\,}r@{\,\,\,\,\,}r@{\,\,}|@{\,\,}r@{\,\,\,\,\,}r@{\,\,}|
 @{\,\,}r@{\,\,\,\,\,}r@{\,\,}|@{\,\,}r@{\,\,\,\,\,}r@{\,\,}|@{\,\,}r@{\,\,\,\,\,}r@{\,\,}|@{\,\,}r@{\,\,\,\,\,}r@{}}
  \hline
$q$ & $\overline{t}_{2}$ &
$q$ & $\overline{t}_{2}$ &
$q$ & $\overline{t}_{2}$ &
$q$ & $\overline{t}_{2}$ &
$q$ & $\overline{t}_{2}$ &
$q$ & $\overline{t}_{2}$ &
$q$ & $\overline{t}_{2}$ &
$q$ & $\overline{t}_{2}^{ \vphantom{H^{L}}}$ \\
 \hline
 2 & 4  &  3 & 4  &  4 & 6  &  5 & 6  &  7 & 6  &  8 & 6  &  9 & 6  &  11 & 7  \\
 13 & 8  &  16 & 9  &  17 & 10  &  19 & 10  &  23 & 10  &  25 & 12  &  27 & 12  &  29 & 13  \\
 31 & 14  &  32 & 14  &  37 & 15  &  41 & 16  &  43 & 16  &  47 & 18  &  49 & 18  &  53 & 18  \\
 59 & 20  &  61 & 20  &  64 & 22  &  67 & 22  &  71 & 22  &  73 & 24  &  79 & 26  &  81 & 26  \\
 83 & 27  &  89 & 28  &  97 & 28  &  101 & 30  &  103 & 31  &  107 & 32  &  109 & 30  &  113 & 33  \\
 121 & 30  &  125 & 30  &  127 & 35  &  128 & 34  &  131 & 36  &  137 & 36  &  139 & 36  &  149 & 39  \\
 151 & 36  &  157 & 40  &  163 & 41  &  167 & 42  &  169 & 42  &  173 & 43  &  179 & 44  &  181 & 44  \\
\hline
 \end{tabular}

 }\newpage

\textbf{Table 4.} Continue 1 \medskip\\\small{
 \renewcommand{\arraystretch}{0.99}
 \begin{tabular}
 {@{}c@{\,\,\,\,\,}c@{\,\,}|@{\,\,}c@{\,\,\,\,\,}c@{\,\,}|@{\,\,}c@{\,\,\,\,\,}c@{\,\,}|@{\,\,}c@{\,\,\,\,\,}c@{\,\,}|@{\,\,}c@{\,\,\,\,\,}c@{\,\,}|@{\,\,}c@{\,\,\,\,\,}c@{\,\,}|@{\,\,}c@{\,\,\,\,\,}c@{\,\,}|@{\,\,}c@{\,\,\,\,\,}c@{}}
  \hline
$q$ & $\overline{t}_{2}$ &
$q$ & $\overline{t}_{2}$ &
$q$ & $\overline{t}_{2}$ &
$q$ & $\overline{t}_{2}$ &
$q$ & $\overline{t}_{2}$ &
$q$ & $\overline{t}_{2}$ &
$q$ & $\overline{t}_{2}$ &
$q$ & $\overline{t}_{2}^{ \vphantom{H^{L}}}$ \\
 \hline
 1901 & 190  &  1907 & 190  &  1913 & 191  &  1931 & 192  &  1933 & 192  &  1949 & 193  &  1951 & 193  &  1973 & 194  \\
 1979 & 195  &  1987 & 195  &  1993 & 195  &  1997 & 196  &  1999 & 196  &  2003 & 196  &  2011 & 197  &  2017 & 197  \\
 2027 & 197  &  2029 & 198  &  2039 & 198  &  2048 & 199  &  2053 & 199  &  2063 & 199  &  2069 & 200  &  2081 & 201  \\
 2083 & 201  &  2087 & 201  &  2089 & 201  &  2099 & 202  &  2111 & 202  &  2113 & 202  &  2129 & 203  &  2131 & 203  \\
 2137 & 204  &  2141 & 204  &  2143 & 204  &  2153 & 205  &  2161 & 205  &  2179 & 206  &  2187 & 207  &  2197 & 207  \\
 2203 & 208  &  2207 & 208  &  2209 & 208  &  2213 & 208  &  2221 & 208  &  2237 & 209  &  2239 & 210  &  2243 & 210  \\
 2251 & 210  &  2267 & 211  &  2269 & 211  &  2273 & 211  &  2281 & 212  &  2287 & 212  &  2293 & 213  &  2297 & 213  \\
\hline
 \end{tabular}

 }\newpage

\textbf{Table 4.} Continue 2 \medskip\\\small{
 \renewcommand{\arraystretch}{0.99}
 \begin{tabular}
 {@{}c@{\,\,\,\,\,}c@{\,\,}|@{\,\,}c@{\,\,\,\,\,}c@{\,\,}|@{\,\,}c@{\,\,\,\,\,}c@{\,\,}|@{\,\,}c@{\,\,\,\,\,}c@{\,\,}|@{\,\,}c@{\,\,\,\,\,}c@{\,\,}|@{\,\,}c@{\,\,\,\,\,}c@{\,\,}|@{\,\,}c@{\,\,\,\,\,}c@{\,\,}|@{\,\,}c@{\,\,\,\,\,}c@{}}
  \hline
$q$ & $\overline{t}_{2}$ &
$q$ & $\overline{t}_{2}$ &
$q$ & $\overline{t}_{2}$ &
$q$ & $\overline{t}_{2}$ &
$q$ & $\overline{t}_{2}$ &
$q$ & $\overline{t}_{2}$ &
$q$ & $\overline{t}_{2}$ &
$q$ & $\overline{t}_{2}^{ \vphantom{H^{L}}}$ \\
 \hline
 4421 & 313  &  4423 & 312  &  4441 & 314  &  4447 & 314  &  4451 & 315  &  4457 & 315  &  4463 & 315  &  4481 & 315  \\
 4483 & 316  &  4489 & 316  &  4493 & 315  &  4507 & 317  &  4513 & 317  &  4517 & 317  &  4519 & 317  &  4523 & 318  \\
 4547 & 319  &  4549 & 319  &  4561 & 319  &  4567 & 319  &  4583 & 320  &  4591 & 320  &  4597 & 321  &  4603 & 321  \\
 4621 & 321  &  4637 & 322  &  4639 & 322  &  4643 & 323  &  4649 & 323  &  4651 & 323  &  4657 & 323  &  4663 & 323  \\
 4673 & 323  &  4679 & 324  &  4691 & 324  &  4703 & 325  &  4721 & 326  &  4723 & 326  &  4729 & 325  &  4733 & 326  \\
 4751 & 327  &  4759 & 327  &  4783 & 328  &  4787 & 328  &  4789 & 327  &  4793 & 329  &  4799 & 329  &  4801 & 329  \\
 4813 & 329  &  4817 & 329  &  4831 & 329  &  4861 & 331  &  4871 & 332  &  4877 & 332  &  4889 & 331  &  4903 & 333  \\
\hline
 \end{tabular}

 }\newpage

\textbf{Table 4.} Continue 3 \medskip\\\small{
 \renewcommand{\arraystretch}{0.99}
 \begin{tabular}
 {@{}c@{\,\,\,\,\,}c@{\,\,}|@{\,\,}c@{\,\,\,\,\,}c@{\,\,}|@{\,\,}c@{\,\,\,\,\,}c@{\,\,}|@{\,\,}c@{\,\,\,\,\,}c@{\,\,}|@{\,\,}c@{\,\,\,\,\,}c@{\,\,}|@{\,\,}c@{\,\,\,\,\,}c@{\,\,}|@{\,\,}c@{\,\,\,\,\,}c@{\,\,}|@{\,\,}c@{\,\,\,\,\,}c@{}}
  \hline
$q$ & $\overline{t}_{2}$ &
$q$ & $\overline{t}_{2}$ &
$q$ & $\overline{t}_{2}$ &
$q$ & $\overline{t}_{2}$ &
$q$ & $\overline{t}_{2}$ &
$q$ & $\overline{t}_{2}$ &
$q$ & $\overline{t}_{2}$ &
$q$ & $\overline{t}_{2}^{ \vphantom{H^{L}}}$ \\
 \hline
 7127 & 414  &  7129 & 414  &  7151 & 415  &  7159 & 415  &  7177 & 416  &  7187 & 415  &  7193 & 415  &  7207 & 417  \\
 7211 & 417  &  7213 & 417  &  7219 & 417  &  7229 & 418  &  7237 & 418  &  7243 & 418  &  7247 & 416  &  7253 & 418  \\
 7283 & 418  &  7297 & 420  &  7307 & 418  &  7309 & 420  &  7321 & 421  &  7331 & 421  &  7333 & 421  &  7349 & 422  \\
 7351 & 422  &  7369 & 422  &  7393 & 423  &  7411 & 423  &  7417 & 423  &  7433 & 424  &  7451 & 422  &  7457 & 424  \\
 7459 & 425  &  7477 & 426  &  7481 & 426  &  7487 & 426  &  7489 & 426  &  7499 & 426  &  7507 & 427  &  7517 & 427  \\
 7523 & 426  &  7529 & 426  &  7537 & 428  &  7541 & 428  &  7547 & 428  &  7549 & 428  &  7559 & 428  &  7561 & 429  \\
 7573 & 429  &  7577 & 428  &  7583 & 429  &  7589 & 428  &  7591 & 430  &  7603 & 430  &  7607 & 430  &  7621 & 431  \\
\hline
 \end{tabular}

 }\newpage

\textbf{Table 5.} The smallest known sizes
 $\overline{t}_{2}=\overline{t}_{2}(2,q)$
 of complete arcs in planes
 $\mathrm{PG}(2,q),$\\ $10001\leq q\leq 100000$, $q$ power
 prime\medskip\\\small{
 \renewcommand{\arraystretch}{0.98}
  \begin{tabular}
 {@{}c@{\,\,\,\,\,}c@{\,\,}|@{\,\,}c@{\,\,\,\,\,}c@{\,\,}|@{\,\,}c@{\,\,\,\,\,}c@{\,\,}|@{\,\,}c@{\,\,\,\,\,}c@{\,\,}|
 @{\,\,}c@{\,\,\,\,\,}c@{\,\,}|@{\,\,}c@{\,\,\,\,\,}c@{\,\,}|@{\,\,}c@{\,\,\,\,\,}c@{}}
  \hline
$q$ & $\overline{t}_{2}$ &
$q$ & $\overline{t}_{2}$ &
$q$ & $\overline{t}_{2}$ &
$q$ & $\overline{t}_{2}$ &
$q$ & $\overline{t}_{2}$ &
$q$ & $\overline{t}_{2}$ &
$q$ & $\overline{t}_{2}^{ \vphantom{H^{L}}}$ \\
 \hline
 10007 & 502  &  10009 & 502  &  10037 & 503  &  10039 & 503  &  10061 & 504  &  10067 & 504  &  10069 & 504  \\
 10079 & 504  &  10091 & 505  &  10093 & 504  &  10099 & 505  &  10103 & 504  &  10111 & 505  &  10133 & 506  \\
 10139 & 506  &  10141 & 506  &  10151 & 504  &  10159 & 505  &  10163 & 506  &  10169 & 507  &  10177 & 507  \\
 10181 & 507  &  10193 & 507  &  10201 & 508  &  10211 & 507  &  10223 & 508  &  10243 & 508  &  10247 & 509  \\
 10253 & 509  &  10259 & 509  &  10267 & 510  &  10271 & 510  &  10273 & 510  &  10289 & 509  &  10301 & 510  \\
 10303 & 510  &  10313 & 511  &  10321 & 511  &  10331 & 511  &  10333 & 511  &  10337 & 512  &  10343 & 511  \\
 10357 & 512  &  10369 & 513  &  10391 & 512  &  10399 & 512  &  10427 & 514  &  10429 & 513  &  10433 & 514  \\
\hline
 \end{tabular}

 }\newpage

\textbf{Table 5.} Continue 1 \medskip\\\small{
 \renewcommand{\arraystretch}{1.00}
 \begin{tabular}
 {@{}c@{\,\,\,\,\,}c@{\,\,}|@{\,\,}c@{\,\,\,\,\,}c@{\,\,}|@{\,\,}c@{\,\,\,\,\,}c@{\,\,}|@{\,\,}c@{\,\,\,\,\,}c@{\,\,}|@{\,\,}c@{\,\,\,\,\,}c@{\,\,}|@{\,\,}c@{\,\,\,\,\,}c@{\,\,}|@{\,\,}c@{\,\,\,\,\,}c@{}}
  \hline
$q$ & $\overline{t}_{2}$ &
$q$ & $\overline{t}_{2}$ &
$q$ & $\overline{t}_{2}$ &
$q$ & $\overline{t}_{2}$ &
$q$ & $\overline{t}_{2}$ &
$q$ & $\overline{t}_{2}$ &
$q$ & $\overline{t}_{2}^{ \vphantom{H^{L}}}$ \\
 \hline
 12601 & 573  &  12611 & 572  &  12613 & 573  &  12619 & 574  &  12637 & 571  &  12641 & 573  &  12647 & 574  \\
 12653 & 575  &  12659 & 573  &  12671 & 575  &  12689 & 575  &  12697 & 576  &  12703 & 576  &  12713 & 576  \\
 12721 & 576  &  12739 & 577  &  12743 & 577  &  12757 & 578  &  12763 & 576  &  12769 & 577  &  12781 & 577  \\
 12791 & 577  &  12799 & 577  &  12809 & 578  &  12821 & 578  &  12823 & 578  &  12829 & 578  &  12841 & 579  \\
 12853 & 578  &  12889 & 581  &  12893 & 580  &  12899 & 581  &  12907 & 581  &  12911 & 580  &  12917 & 580  \\
 12919 & 580  &  12923 & 581  &  12941 & 582  &  12953 & 582  &  12959 & 581  &  12967 & 582  &  12973 & 582  \\
 12979 & 582  &  12983 & 583  &  13001 & 582  &  13003 & 581  &  13007 & 583  &  13009 & 583  &  13033 & 584  \\
\hline
 \end{tabular}

 }\newpage

\textbf{Table 5.} Continue 2 \medskip\\\small{
 \renewcommand{\arraystretch}{1.00}
 \begin{tabular}
 {@{}c@{\,\,\,\,\,}c@{\,\,}|@{\,\,}c@{\,\,\,\,\,}c@{\,\,}|@{\,\,}c@{\,\,\,\,\,}c@{\,\,}|@{\,\,}c@{\,\,\,\,\,}c@{\,\,}|@{\,\,}c@{\,\,\,\,\,}c@{\,\,}|@{\,\,}c@{\,\,\,\,\,}c@{\,\,}|@{\,\,}c@{\,\,\,\,\,}c@{}}
  \hline
$q$ & $\overline{t}_{2}$ &
$q$ & $\overline{t}_{2}$ &
$q$ & $\overline{t}_{2}$ &
$q$ & $\overline{t}_{2}$ &
$q$ & $\overline{t}_{2}$ &
$q$ & $\overline{t}_{2}$ &
$q$ & $\overline{t}_{2}^{ \vphantom{H^{L}}}$ \\
 \hline
 15271 & 638  &  15277 & 639  &  15287 & 639  &  15289 & 641  &  15299 & 639  &  15307 & 639  &  15313 & 642  \\
 15319 & 639  &  15329 & 640  &  15331 & 638  &  15349 & 639  &  15359 & 641  &  15361 & 641  &  15373 & 641  \\
 15377 & 643  &  15383 & 641  &  15391 & 639  &  15401 & 641  &  15413 & 641  &  15427 & 641  &  15439 & 643  \\
 15443 & 644  &  15451 & 644  &  15461 & 643  &  15467 & 643  &  15473 & 644  &  15493 & 643  &  15497 & 644  \\
 15511 & 643  &  15527 & 644  &  15541 & 645  &  15551 & 645  &  15559 & 647  &  15569 & 645  &  15581 & 646  \\
 15583 & 645  &  15601 & 645  &  15607 & 647  &  15619 & 647  &  15625 & 647  &  15629 & 647  &  15641 & 645  \\
 15643 & 648  &  15647 & 648  &  15649 & 646  &  15661 & 649  &  15667 & 648  &  15671 & 648  &  15679 & 648  \\
\hline
 \end{tabular}

 }\newpage

\textbf{Table 5.} Continue 3 \medskip\\\small{
 \renewcommand{\arraystretch}{1.00}
 \begin{tabular}
 {@{}c@{\,\,\,\,\,}c@{\,\,}|@{\,\,}c@{\,\,\,\,\,}c@{\,\,}|@{\,\,}c@{\,\,\,\,\,}c@{\,\,}|@{\,\,}c@{\,\,\,\,\,}c@{\,\,}|@{\,\,}c@{\,\,\,\,\,}c@{\,\,}|@{\,\,}c@{\,\,\,\,\,}c@{\,\,}|@{\,\,}c@{\,\,\,\,\,}c@{}}
  \hline
$q$ & $\overline{t}_{2}$ &
$q$ & $\overline{t}_{2}$ &
$q$ & $\overline{t}_{2}$ &
$q$ & $\overline{t}_{2}$ &
$q$ & $\overline{t}_{2}$ &
$q$ & $\overline{t}_{2}$ &
$q$ & $\overline{t}_{2}^{ \vphantom{H^{L}}}$ \\
 \hline
 17957 & 700  &  17959 & 701  &  17971 & 701  &  17977 & 700  &  17981 & 701  &  17987 & 701  &  17989 & 701  \\
 18013 & 700  &  18041 & 702  &  18043 & 703  &  18047 & 702  &  18049 & 702  &  18059 & 702  &  18061 & 702  \\
 18077 & 702  &  18089 & 703  &  18097 & 702  &  18119 & 702  &  18121 & 704  &  18127 & 703  &  18131 & 704  \\
 18133 & 704  &  18143 & 705  &  18149 & 704  &  18169 & 705  &  18181 & 706  &  18191 & 705  &  18199 & 705  \\
 18211 & 706  &  18217 & 708  &  18223 & 706  &  18229 & 705  &  18233 & 706  &  18251 & 705  &  18253 & 707  \\
 18257 & 706  &  18269 & 707  &  18287 & 706  &  18289 & 708  &  18301 & 705  &  18307 & 707  &  18311 & 708  \\
 18313 & 708  &  18329 & 709  &  18341 & 710  &  18353 & 710  &  18367 & 707  &  18371 & 709  &  18379 & 709  \\
\hline
 \end{tabular}

 }\newpage

\textbf{Table 5.} Continue 4 \medskip\\\small{
 \renewcommand{\arraystretch}{1.00}
 \begin{tabular}
 {@{}c@{\,\,\,\,\,}c@{\,\,}|@{\,\,}c@{\,\,\,\,\,}c@{\,\,}|@{\,\,}c@{\,\,\,\,\,}c@{\,\,}|@{\,\,}c@{\,\,\,\,\,}c@{\,\,}|@{\,\,}c@{\,\,\,\,\,}c@{\,\,}|@{\,\,}c@{\,\,\,\,\,}c@{\,\,}|@{\,\,}c@{\,\,\,\,\,}c@{}}
  \hline
$q$ & $\overline{t}_{2}$ &
$q$ & $\overline{t}_{2}$ &
$q$ & $\overline{t}_{2}$ &
$q$ & $\overline{t}_{2}$ &
$q$ & $\overline{t}_{2}$ &
$q$ & $\overline{t}_{2}$ &
$q$ & $\overline{t}_{2}^{ \vphantom{H^{L}}}$ \\
 \hline
 20743 & 760  &  20747 & 761  &  20749 & 760  &  20753 & 759  &  20759 & 761  &  20771 & 760  &  20773 & 760  \\
 20789 & 761  &  20807 & 759  &  20809 & 762  &  20849 & 762  &  20857 & 763  &  20873 & 763  &  20879 & 761  \\
 20887 & 764  &  20897 & 763  &  20899 & 765  &  20903 & 762  &  20921 & 762  &  20929 & 764  &  20939 & 763  \\
 20947 & 764  &  20959 & 763  &  20963 & 764  &  20981 & 762  &  20983 & 764  &  21001 & 765  &  21011 & 766  \\
 21013 & 765  &  21017 & 765  &  21019 & 765  &  21023 & 764  &  21031 & 766  &  21059 & 766  &  21061 & 766  \\
 21067 & 767  &  21089 & 768  &  21101 & 765  &  21107 & 768  &  21121 & 768  &  21139 & 766  &  21143 & 769  \\
 21149 & 766  &  21157 & 768  &  21163 & 768  &  21169 & 766  &  21179 & 768  &  21187 & 768  &  21191 & 770  \\
\hline
 \end{tabular}

 }\newpage

\textbf{Table 5.} Continue 5 \medskip\\\small{
 \renewcommand{\arraystretch}{1.00}
 \begin{tabular}
 {@{}c@{\,\,\,\,\,}c@{\,\,}|@{\,\,}c@{\,\,\,\,\,}c@{\,\,}|@{\,\,}c@{\,\,\,\,\,}c@{\,\,}|@{\,\,}c@{\,\,\,\,\,}c@{\,\,}|@{\,\,}c@{\,\,\,\,\,}c@{\,\,}|@{\,\,}c@{\,\,\,\,\,}c@{\,\,}|@{\,\,}c@{\,\,\,\,\,}c@{}}
  \hline
$q$ & $\overline{t}_{2}$ &
$q$ & $\overline{t}_{2}$ &
$q$ & $\overline{t}_{2}$ &
$q$ & $\overline{t}_{2}$ &
$q$ & $\overline{t}_{2}$ &
$q$ & $\overline{t}_{2}$ &
$q$ & $\overline{t}_{2}^{ \vphantom{H^{L}}}$ \\
 \hline
 23497 & 815  &  23509 & 814  &  23531 & 817  &  23537 & 815  &  23539 & 816  &  23549 & 814  &  23557 & 817  \\
 23561 & 813  &  23563 & 817  &  23567 & 817  &  23581 & 817  &  23593 & 816  &  23599 & 817  &  23603 & 817  \\
 23609 & 817  &  23623 & 818  &  23627 & 816  &  23629 & 818  &  23633 & 816  &  23663 & 819  &  23669 & 818  \\
 23671 & 819  &  23677 & 820  &  23687 & 817  &  23689 & 817  &  23719 & 819  &  23741 & 819  &  23743 & 821  \\
 23747 & 819  &  23753 & 818  &  23761 & 822  &  23767 & 823  &  23773 & 819  &  23789 & 821  &  23801 & 821  \\
 23813 & 821  &  23819 & 819  &  23827 & 822  &  23831 & 820  &  23833 & 820  &  23857 & 824  &  23869 & 822  \\
 23873 & 822  &  23879 & 824  &  23887 & 823  &  23893 & 823  &  23899 & 823  &  23909 & 822  &  23911 & 823  \\
\hline
 \end{tabular}

 }\newpage

\textbf{Table 5.} Continue 6 \medskip\\\small{
 \renewcommand{\arraystretch}{1.00}
 \begin{tabular}
 {@{}c@{\,\,\,\,\,}c@{\,\,}|@{\,\,}c@{\,\,\,\,\,}c@{\,\,}|@{\,\,}c@{\,\,\,\,\,}c@{\,\,}|@{\,\,}c@{\,\,\,\,\,}c@{\,\,}|@{\,\,}c@{\,\,\,\,\,}c@{\,\,}|@{\,\,}c@{\,\,\,\,\,}c@{\,\,}|@{\,\,}c@{\,\,\,\,\,}c@{}}
  \hline
$q$ & $\overline{t}_{2}$ &
$q$ & $\overline{t}_{2}$ &
$q$ & $\overline{t}_{2}$ &
$q$ & $\overline{t}_{2}$ &
$q$ & $\overline{t}_{2}$ &
$q$ & $\overline{t}_{2}$ &
$q$ & $\overline{t}_{2}^{ \vphantom{H^{L}}}$ \\
 \hline
 26309 & 869  &  26317 & 869  &  26321 & 871  &  26339 & 871  &  26347 & 870  &  26357 & 870  &  26371 & 868  \\
 26387 & 872  &  26393 & 871  &  26399 & 870  &  26407 & 871  &  26417 & 869  &  26423 & 869  &  26431 & 871  \\
 26437 & 872  &  26449 & 870  &  26459 & 871  &  26479 & 872  &  26489 & 869  &  26497 & 874  &  26501 & 872  \\
 26513 & 874  &  26539 & 873  &  26557 & 873  &  26561 & 872  &  26569 & 873  &  26573 & 875  &  26591 & 874  \\
 26597 & 876  &  26627 & 871  &  26633 & 876  &  26641 & 874  &  26647 & 875  &  26669 & 876  &  26681 & 875  \\
 26683 & 875  &  26687 & 878  &  26693 & 875  &  26699 & 874  &  26701 & 876  &  26711 & 876  &  26713 & 877  \\
 26717 & 876  &  26723 & 876  &  26729 & 877  &  26731 & 875  &  26737 & 875  &  26759 & 880  &  26777 & 878  \\
\hline
 \end{tabular}

 }\newpage

\textbf{Table 5.} Continue 7 \medskip\\\small{
 \renewcommand{\arraystretch}{1.00}
 \begin{tabular}
 {@{}c@{\,\,\,\,\,}c@{\,\,}|@{\,\,}c@{\,\,\,\,\,}c@{\,\,}|@{\,\,}c@{\,\,\,\,\,}c@{\,\,}|@{\,\,}c@{\,\,\,\,\,}c@{\,\,}|@{\,\,}c@{\,\,\,\,\,}c@{\,\,}|@{\,\,}c@{\,\,\,\,\,}c@{\,\,}|@{\,\,}c@{\,\,\,\,\,}c@{}}
  \hline
$q$ & $\overline{t}_{2}$ &
$q$ & $\overline{t}_{2}$ &
$q$ & $\overline{t}_{2}$ &
$q$ & $\overline{t}_{2}$ &
$q$ & $\overline{t}_{2}$ &
$q$ & $\overline{t}_{2}$ &
$q$ & $\overline{t}_{2}^{ \vphantom{H^{L}}}$ \\
 \hline
 29147 & 924  &  29153 & 922  &  29167 & 919  &  29173 & 919  &  29179 & 922  &  29191 & 923  &  29201 & 923  \\
 29207 & 920  &  29209 & 920  &  29221 & 919  &  29231 & 921  &  29243 & 924  &  29251 & 921  &  29269 & 921  \\
 29287 & 925  &  29297 & 923  &  29303 & 924  &  29311 & 924  &  29327 & 923  &  29333 & 922  &  29339 & 925  \\
 29347 & 924  &  29363 & 925  &  29383 & 923  &  29387 & 924  &  29389 & 926  &  29399 & 922  &  29401 & 923  \\
 29411 & 924  &  29423 & 927  &  29429 & 928  &  29437 & 926  &  29443 & 924  &  29453 & 927  &  29473 & 928  \\
 29483 & 928  &  29501 & 926  &  29527 & 928  &  29531 & 927  &  29537 & 929  &  29567 & 928  &  29569 & 926  \\
 29573 & 928  &  29581 & 928  &  29587 & 927  &  29599 & 929  &  29611 & 927  &  29629 & 929  &  29633 & 928  \\
\hline
 \end{tabular}

 }\newpage

\textbf{Table 5.} Continue 8 \medskip\\\small{
 \renewcommand{\arraystretch}{1.00}
 \begin{tabular}
 {@{}c@{\,\,\,\,\,}r@{\,\,}|@{\,\,}c@{\,\,\,\,\,}r@{\,\,}|@{\,\,}c@{\,\,\,\,\,}r@{\,\,}|@{\,\,}c@{\,\,\,\,\,}r@{\,\,}|
 @{\,\,}c@{\,\,\,\,\,}r@{\,\,}|@{\,\,}c@{\,\,\,\,\,}r@{\,\,}|@{\,\,}c@{\,\,\,\,\,}r@{}}
  \hline
$q$ & $\overline{t}_{2}$ &
$q$ & $\overline{t}_{2}$ &
$q$ & $\overline{t}_{2}$ &
$q$ & $\overline{t}_{2}$ &
$q$ & $\overline{t}_{2}$ &
$q$ & $\overline{t}_{2}$ &
$q$ & $\overline{t}_{2}^{ \vphantom{H^{L}}}$ \\
 \hline
 32099 & 972  &  32117 & 972  &  32119 & 972  &  32141 & 973  &  32143 & 972  &  32159 & 974  &  32173 & 974  \\
 32183 & 972  &  32189 & 972  &  32191 & 975  &  32203 & 973  &  32213 & 972  &  32233 & 975  &  32237 & 977  \\
 32251 & 978  &  32257 & 974  &  32261 & 977  &  32297 & 975  &  32299 & 974  &  32303 & 977  &  32309 & 974  \\
 32321 & 975  &  32323 & 977  &  32327 & 975  &  32341 & 976  &  32353 & 974  &  32359 & 976  &  32363 & 978  \\
 32369 & 977  &  32371 & 975  &  32377 & 977  &  32381 & 978  &  32401 & 979  &  32411 & 976  &  32413 & 978  \\
 32423 & 977  &  32429 & 977  &  32441 & 978  &  32443 & 977  &  32467 & 980  &  32479 & 979  &  32491 & 981  \\
 32497 & 982  &  32503 & 979  &  32507 & 981  &  32531 & 981  &  32533 & 977  &  32537 & 980  &  32561 & 978  \\
\hline
 \end{tabular}

 }\newpage

\textbf{Table 5.} Continue 9 \medskip\\\small{
 \renewcommand{\arraystretch}{1.00}
 \begin{tabular}
 {@{}c@{\,\,\,\,\,}c@{\,\,}|@{\,\,}c@{\,\,\,\,\,}c@{\,\,}|@{\,\,}c@{\,\,\,\,\,}c@{\,\,}|@{\,\,}c@{\,\,\,\,\,}c@{\,\,}|@{\,\,}c@{\,\,\,\,\,}c@{\,\,}|@{\,\,}c@{\,\,\,\,\,}c@{\,\,}|@{\,\,}c@{\,\,\,\,\,}c@{}}
  \hline
$q$ & $\overline{t}_{2}$ &
$q$ & $\overline{t}_{2}$ &
$q$ & $\overline{t}_{2}$ &
$q$ & $\overline{t}_{2}$ &
$q$ & $\overline{t}_{2}$ &
$q$ & $\overline{t}_{2}$ &
$q$ & $\overline{t}_{2}^{ \vphantom{H^{L}}}$ \\
 \hline
 34877 & 1019  &  34883 & 1018  &  34897 & 1018  &  34913 & 1020  &  34919 & 1020  &  34939 & 1018  &  34949 & 1021  \\
 34961 & 1022  &  34963 & 1020  &  34981 & 1022  &  35023 & 1019  &  35027 & 1018  &  35051 & 1022  &  35053 & 1019  \\
 35059 & 1021  &  35069 & 1021  &  35081 & 1020  &  35083 & 1024  &  35089 & 1021  &  35099 & 1022  &  35107 & 1022  \\
 35111 & 1023  &  35117 & 1021  &  35129 & 1020  &  35141 & 1022  &  35149 & 1022  &  35153 & 1023  &  35159 & 1025  \\
 35171 & 1022  &  35201 & 1023  &  35221 & 1027  &  35227 & 1024  &  35251 & 1027  &  35257 & 1025  &  35267 & 1026  \\
 35279 & 1024  &  35281 & 1024  &  35291 & 1024  &  35311 & 1022  &  35317 & 1028  &  35323 & 1027  &  35327 & 1026  \\
 35339 & 1023  &  35353 & 1026  &  35363 & 1025  &  35381 & 1028  &  35393 & 1028  &  35401 & 1028  &  35407 & 1028  \\
\hline
 \end{tabular}

 }\newpage

\textbf{Table 5.} Continue 10 \medskip\\\small{
 \renewcommand{\arraystretch}{1.00}
 \begin{tabular}
 {@{}c@{\,\,\,\,\,}c@{\,\,}|@{\,\,}c@{\,\,\,\,\,}c@{\,\,}|@{\,\,}c@{\,\,\,\,\,}c@{\,\,}|@{\,\,}c@{\,\,\,\,\,}c@{\,\,}|@{\,\,}c@{\,\,\,\,\,}c@{\,\,}|@{\,\,}c@{\,\,\,\,\,}c@{\,\,}|@{\,\,}c@{\,\,\,\,\,}c@{}}
  \hline
$q$ & $\overline{t}_{2}$ &
$q$ & $\overline{t}_{2}$ &
$q$ & $\overline{t}_{2}$ &
$q$ & $\overline{t}_{2}$ &
$q$ & $\overline{t}_{2}$ &
$q$ & $\overline{t}_{2}$ &
$q$ & $\overline{t}_{2}^{ \vphantom{H^{L}}}$ \\
 \hline
 37831 & 1065  &  37847 & 1065  &  37853 & 1064  &  37861 & 1064  &  37871 & 1067  &  37879 & 1064  &  37889 & 1067  \\
 37897 & 1065  &  37907 & 1069  &  37951 & 1069  &  37957 & 1067  &  37963 & 1072  &  37967 & 1062  &  37987 & 1066  \\
 37991 & 1068  &  37993 & 1069  &  37997 & 1070  &  38011 & 1070  &  38039 & 1073  &  38047 & 1072  &  38053 & 1072  \\
 38069 & 1068  &  38083 & 1073  &  38113 & 1072  &  38119 & 1072  &  38149 & 1070  &  38153 & 1070  &  38167 & 1073  \\
 38177 & 1071  &  38183 & 1071  &  38189 & 1072  &  38197 & 1077  &  38201 & 1071  &  38219 & 1069  &  38231 & 1070  \\
 38237 & 1074  &  38239 & 1072  &  38261 & 1072  &  38273 & 1072  &  38281 & 1073  &  38287 & 1072  &  38299 & 1075  \\
 38303 & 1073  &  38317 & 1076  &  38321 & 1072  &  38327 & 1073  &  38329 & 1071  &  38333 & 1073  &  38351 & 1075  \\
\hline
 \end{tabular}

 }\newpage

\textbf{Table 5.} Continue 11 \medskip\\\small{
 \renewcommand{\arraystretch}{1.00}
 \begin{tabular}
 {@{}c@{\,\,\,\,\,}c@{\,\,}|@{\,\,}c@{\,\,\,\,\,}c@{\,\,}|@{\,\,}c@{\,\,\,\,\,}c@{\,\,}|@{\,\,}c@{\,\,\,\,\,}c@{\,\,}|@{\,\,}c@{\,\,\,\,\,}c@{\,\,}|@{\,\,}c@{\,\,\,\,\,}c@{\,\,}|@{\,\,}c@{\,\,\,\,\,}c@{}}
  \hline
$q$ & $\overline{t}_{2}$ &
$q$ & $\overline{t}_{2}$ &
$q$ & $\overline{t}_{2}$ &
$q$ & $\overline{t}_{2}$ &
$q$ & $\overline{t}_{2}$ &
$q$ & $\overline{t}_{2}$ &
$q$ & $\overline{t}_{2}^{ \vphantom{H^{L}}}$ \\
 \hline
 40853 & 1113  &  40867 & 1118  &  40879 & 1116  &  40883 & 1115  &  40897 & 1113  &  40903 & 1115  &  40927 & 1118  \\
 40933 & 1115  &  40939 & 1114  &  40949 & 1115  &  40961 & 1117  &  40973 & 1117  &  40993 & 1116  &  41011 & 1115  \\
 41017 & 1115  &  41023 & 1115  &  41039 & 1114  &  41047 & 1119  &  41051 & 1116  &  41057 & 1118  &  41077 & 1119  \\
 41081 & 1117  &  41113 & 1121  &  41117 & 1114  &  41131 & 1115  &  41141 & 1120  &  41143 & 1120  &  41149 & 1119  \\
 41161 & 1120  &  41177 & 1119  &  41179 & 1119  &  41183 & 1123  &  41189 & 1119  &  41201 & 1122  &  41203 & 1117  \\
 41213 & 1120  &  41221 & 1119  &  41227 & 1121  &  41231 & 1120  &  41233 & 1118  &  41243 & 1123  &  41257 & 1118  \\
 41263 & 1124  &  41269 & 1118  &  41281 & 1122  &  41299 & 1120  &  41333 & 1121  &  41341 & 1120  &  41351 & 1126  \\
\hline
 \end{tabular}

 }\newpage

\textbf{Table 5.} Continue 12 \medskip\\\small{
 \renewcommand{\arraystretch}{1.00}
 \begin{tabular}
 {@{}c@{\,\,\,\,\,}c@{\,\,}|@{\,\,}c@{\,\,\,\,\,}c@{\,\,}|@{\,\,}c@{\,\,\,\,\,}c@{\,\,}|@{\,\,}c@{\,\,\,\,\,}c@{\,\,}|@{\,\,}c@{\,\,\,\,\,}c@{\,\,}|@{\,\,}c@{\,\,\,\,\,}c@{\,\,}|@{\,\,}c@{\,\,\,\,\,}c@{}}
  \hline
$q$ & $\overline{t}_{2}$ &
$q$ & $\overline{t}_{2}$ &
$q$ & $\overline{t}_{2}$ &
$q$ & $\overline{t}_{2}$ &
$q$ & $\overline{t}_{2}$ &
$q$ & $\overline{t}_{2}$ &
$q$ & $\overline{t}_{2}^{ \vphantom{H^{L}}}$ \\
 \hline
 43781 & 1160  &  43783 & 1159  &  43787 & 1156  &  43789 & 1157  &  43793 & 1158  &  43801 & 1161  &  43853 & 1160  \\
 43867 & 1165  &  43889 & 1164  &  43891 & 1160  &  43913 & 1163  &  43933 & 1163  &  43943 & 1161  &  43951 & 1161  \\
 43961 & 1162  &  43963 & 1162  &  43969 & 1161  &  43973 & 1164  &  43987 & 1162  &  43991 & 1162  &  43997 & 1165  \\
 44017 & 1162  &  44021 & 1161  &  44027 & 1163  &  44029 & 1161  &  44041 & 1163  &  44053 & 1161  &  44059 & 1163  \\
 44071 & 1166  &  44087 & 1162  &  44089 & 1165  &  44101 & 1164  &  44111 & 1163  &  44119 & 1160  &  44123 & 1163  \\
 44129 & 1162  &  44131 & 1166  &  44159 & 1161  &  44171 & 1162  &  44179 & 1163  &  44189 & 1162  &  44201 & 1162  \\
 44203 & 1169  &  44207 & 1164  &  44221 & 1161  &  44249 & 1163  &  44257 & 1164  &  44263 & 1167  &  44267 & 1165  \\
\hline
 \end{tabular}

 }\newpage

\textbf{Table 5.} Continue 13 \medskip\\\small{
 \renewcommand{\arraystretch}{1.00}
 \begin{tabular}
 {@{}c@{\,\,\,\,\,}c@{\,\,}|@{\,\,}c@{\,\,\,\,\,}c@{\,\,}|@{\,\,}c@{\,\,\,\,\,}c@{\,\,}|@{\,\,}c@{\,\,\,\,\,}c@{\,\,}|@{\,\,}c@{\,\,\,\,\,}c@{\,\,}|@{\,\,}c@{\,\,\,\,\,}c@{\,\,}|@{\,\,}c@{\,\,\,\,\,}c@{}}
  \hline
$q$ & $\overline{t}_{2}$ &
$q$ & $\overline{t}_{2}$ &
$q$ & $\overline{t}_{2}$ &
$q$ & $\overline{t}_{2}$ &
$q$ & $\overline{t}_{2}$ &
$q$ & $\overline{t}_{2}$ &
$q$ & $\overline{t}_{2}^{ \vphantom{H^{L}}}$ \\
 \hline
 46819 & 1204  &  46829 & 1205  &  46831 & 1208  &  46853 & 1205  &  46861 & 1200  &  46867 & 1198  &  46877 & 1204  \\
 46889 & 1200  &  46901 & 1204  &  46919 & 1204  &  46933 & 1205  &  46957 & 1205  &  46993 & 1208  &  46997 & 1204  \\
 47017 & 1205  &  47041 & 1209  &  47051 & 1206  &  47057 & 1208  &  47059 & 1208  &  47087 & 1208  &  47093 & 1210  \\
 47111 & 1206  &  47119 & 1208  &  47123 & 1210  &  47129 & 1206  &  47137 & 1205  &  47143 & 1208  &  47147 & 1210  \\
 47149 & 1209  &  47161 & 1205  &  47189 & 1211  &  47207 & 1210  &  47221 & 1206  &  47237 & 1206  &  47251 & 1208  \\
 47269 & 1208  &  47279 & 1210  &  47287 & 1212  &  47293 & 1211  &  47297 & 1205  &  47303 & 1208  &  47309 & 1213  \\
 47317 & 1209  &  47339 & 1206  &  47351 & 1212  &  47353 & 1206  &  47363 & 1211  &  47381 & 1211  &  47387 & 1205  \\
\hline
 \end{tabular}

 }\newpage

\textbf{Table 5.} Continue 14 \medskip\\\small{
 \renewcommand{\arraystretch}{1.00}
 \begin{tabular}
 {@{}c@{\,\,\,\,\,}c@{\,\,}|@{\,\,}c@{\,\,\,\,\,}c@{\,\,}|@{\,\,}c@{\,\,\,\,\,}c@{\,\,}|@{\,\,}c@{\,\,\,\,\,}c@{\,\,}|@{\,\,}c@{\,\,\,\,\,}c@{\,\,}|@{\,\,}c@{\,\,\,\,\,}c@{\,\,}|@{\,\,}c@{\,\,\,\,\,}c@{}}
  \hline
$q$ & $\overline{t}_{2}$ &
$q$ & $\overline{t}_{2}$ &
$q$ & $\overline{t}_{2}$ &
$q$ & $\overline{t}_{2}$ &
$q$ & $\overline{t}_{2}$ &
$q$ & $\overline{t}_{2}$ &
$q$ & $\overline{t}_{2}^{ \vphantom{H^{L}}}$ \\
 \hline
 49823 & 1244  &  49831 & 1241  &  49843 & 1248  &  49853 & 1250  &  49871 & 1248  &  49877 & 1245  &  49891 & 1248  \\
 49919 & 1244  &  49921 & 1246  &  49927 & 1246  &  49937 & 1247  &  49939 & 1245  &  49943 & 1247  &  49957 & 1248  \\
 49991 & 1245  &  49993 & 1252  &  49999 & 1246  &  50021 & 1247  &  50023 & 1248  &  50033 & 1248  &  50047 & 1249  \\
 50051 & 1252  &  50053 & 1244  &  50069 & 1247  &  50077 & 1253  &  50087 & 1248  &  50093 & 1244  &  50101 & 1250  \\
 50111 & 1253  &  50119 & 1251  &  50123 & 1250  &  50129 & 1253  &  50131 & 1248  &  50147 & 1252  &  50153 & 1254  \\
 50159 & 1247  &  50177 & 1247  &  50207 & 1247  &  50221 & 1254  &  50227 & 1255  &  50231 & 1255  &  50261 & 1253  \\
 50263 & 1251  &  50273 & 1249  &  50287 & 1252  &  50291 & 1253  &  50311 & 1253  &  50321 & 1248  &  50329 & 1251  \\
\hline
 \end{tabular}

 }\newpage

\textbf{Table 5.} Continue 15 \medskip\\\small{
 \renewcommand{\arraystretch}{1.00}
 \begin{tabular}
 {@{}c@{\,\,\,\,\,}c@{\,\,}|@{\,\,}c@{\,\,\,\,\,}c@{\,\,}|@{\,\,}c@{\,\,\,\,\,}c@{\,\,}|@{\,\,}c@{\,\,\,\,\,}c@{\,\,}|@{\,\,}c@{\,\,\,\,\,}c@{\,\,}|@{\,\,}c@{\,\,\,\,\,}c@{\,\,}|@{\,\,}c@{\,\,\,\,\,}c@{}}
  \hline
$q$ & $\overline{t}_{2}$ &
$q$ & $\overline{t}_{2}$ &
$q$ & $\overline{t}_{2}$ &
$q$ & $\overline{t}_{2}$ &
$q$ & $\overline{t}_{2}$ &
$q$ & $\overline{t}_{2}$ &
$q$ & $\overline{t}_{2}^{ \vphantom{H^{L}}}$ \\
 \hline
 52861 & 1288  &  52879 & 1290  &  52883 & 1289  &  52889 & 1288  &  52901 & 1284  &  52903 & 1287  &  52919 & 1285  \\
 52937 & 1287  &  52951 & 1283  &  52957 & 1288  &  52963 & 1288  &  52967 & 1284  &  52973 & 1290  &  52981 & 1289  \\
 52999 & 1284  &  53003 & 1290  &  53017 & 1289  &  53047 & 1291  &  53051 & 1285  &  53069 & 1291  &  53077 & 1289  \\
 53087 & 1295  &  53089 & 1290  &  53093 & 1292  &  53101 & 1296  &  53113 & 1294  &  53117 & 1290  &  53129 & 1300  \\
 53147 & 1292  &  53149 & 1291  &  53161 & 1289  &  53171 & 1295  &  53173 & 1291  &  53189 & 1291  &  53197 & 1294  \\
 53201 & 1294  &  53231 & 1286  &  53233 & 1295  &  53239 & 1297  &  53267 & 1293  &  53269 & 1289  &  53279 & 1291  \\
 53281 & 1293  &  53299 & 1294  &  53309 & 1295  &  53323 & 1298  &  53327 & 1297  &  53353 & 1292  &  53359 & 1292  \\
\hline
 \end{tabular}

 }\newpage

\textbf{Table 5.} Continue 16 \medskip\\\small{
 \renewcommand{\arraystretch}{1.00}
 \begin{tabular}
 {@{}c@{\,\,\,\,\,}c@{\,\,}|@{\,\,}c@{\,\,\,\,\,}c@{\,\,}|@{\,\,}c@{\,\,\,\,\,}c@{\,\,}|@{\,\,}c@{\,\,\,\,\,}c@{\,\,}|@{\,\,}c@{\,\,\,\,\,}c@{\,\,}|@{\,\,}c@{\,\,\,\,\,}c@{\,\,}|@{\,\,}c@{\,\,\,\,\,}c@{}}
  \hline
$q$ & $\overline{t}_{2}$ &
$q$ & $\overline{t}_{2}$ &
$q$ & $\overline{t}_{2}$ &
$q$ & $\overline{t}_{2}$ &
$q$ & $\overline{t}_{2}$ &
$q$ & $\overline{t}_{2}$ &
$q$ & $\overline{t}_{2}^{ \vphantom{H^{L}}}$ \\
 \hline
 55897 & 1327  &  55901 & 1330  &  55903 & 1333  &  55921 & 1331  &  55927 & 1330  &  55931 & 1330  &  55933 & 1326  \\
 55949 & 1330  &  55967 & 1331  &  55987 & 1336  &  55997 & 1337  &  56003 & 1328  &  56009 & 1328  &  56039 & 1325  \\
 56041 & 1333  &  56053 & 1326  &  56081 & 1334  &  56087 & 1331  &  56093 & 1326  &  56099 & 1332  &  56101 & 1330  \\
 56113 & 1328  &  56123 & 1332  &  56131 & 1333  &  56149 & 1330  &  56167 & 1329  &  56171 & 1331  &  56179 & 1333  \\
 56197 & 1329  &  56207 & 1326  &  56209 & 1337  &  56237 & 1333  &  56239 & 1333  &  56249 & 1333  &  56263 & 1330  \\
 56267 & 1334  &  56269 & 1336  &  56299 & 1333  &  56311 & 1335  &  56333 & 1332  &  56359 & 1333  &  56369 & 1335  \\
 56377 & 1339  &  56383 & 1337  &  56393 & 1337  &  56401 & 1336  &  56417 & 1331  &  56431 & 1333  &  56437 & 1330  \\
\hline
 \end{tabular}

 }\newpage

\textbf{Table 5.} Continue 17 \medskip\\\small{
 \renewcommand{\arraystretch}{1.00}
 \begin{tabular}
 {@{}c@{\,\,\,\,\,}c@{\,\,}|@{\,\,}c@{\,\,\,\,\,}c@{\,\,}|@{\,\,}c@{\,\,\,\,\,}c@{\,\,}|@{\,\,}c@{\,\,\,\,\,}c@{\,\,}|@{\,\,}c@{\,\,\,\,\,}c@{\,\,}|@{\,\,}c@{\,\,\,\,\,}c@{\,\,}|@{\,\,}c@{\,\,\,\,\,}c@{}}
  \hline
$q$ & $\overline{t}_{2}$ &
$q$ & $\overline{t}_{2}$ &
$q$ & $\overline{t}_{2}$ &
$q$ & $\overline{t}_{2}$ &
$q$ & $\overline{t}_{2}$ &
$q$ & $\overline{t}_{2}$ &
$q$ & $\overline{t}_{2}^{ \vphantom{H^{L}}}$ \\
 \hline
 58897 & 1369  &  58901 & 1369  &  58907 & 1364  &  58909 & 1369  &  58913 & 1372  &  58921 & 1369  &  58937 & 1368  \\
 58943 & 1373  &  58963 & 1363  &  58967 & 1369  &  58979 & 1366  &  58991 & 1373  &  58997 & 1365  &  59009 & 1374  \\
 59011 & 1372  &  59021 & 1367  &  59023 & 1368  &  59029 & 1372  &  59049 & 1370  &  59051 & 1368  &  59053 & 1370  \\
 59063 & 1371  &  59069 & 1368  &  59077 & 1378  &  59083 & 1364  &  59093 & 1374  &  59107 & 1370  &  59113 & 1371  \\
 59119 & 1371  &  59123 & 1370  &  59141 & 1366  &  59149 & 1369  &  59159 & 1380  &  59167 & 1373  &  59183 & 1372  \\
 59197 & 1368  &  59207 & 1372  &  59209 & 1369  &  59219 & 1369  &  59221 & 1376  &  59233 & 1370  &  59239 & 1369  \\
 59243 & 1374  &  59263 & 1376  &  59273 & 1374  &  59281 & 1376  &  59333 & 1376  &  59341 & 1373  &  59351 & 1377  \\
\hline
 \end{tabular}

 }\newpage

\textbf{Table 5.} Continue 18 \medskip\\\small{
 \renewcommand{\arraystretch}{1.00}
 \begin{tabular}
 {@{}c@{\,\,\,\,\,}c@{\,\,}|@{\,\,}c@{\,\,\,\,\,}c@{\,\,}|@{\,\,}c@{\,\,\,\,\,}c@{\,\,}|@{\,\,}c@{\,\,\,\,\,}c@{\,\,}|@{\,\,}c@{\,\,\,\,\,}c@{\,\,}|@{\,\,}c@{\,\,\,\,\,}c@{\,\,}|@{\,\,}c@{\,\,\,\,\,}c@{}}
  \hline
$q$ & $\overline{t}_{2}$ &
$q$ & $\overline{t}_{2}$ &
$q$ & $\overline{t}_{2}$ &
$q$ & $\overline{t}_{2}$ &
$q$ & $\overline{t}_{2}$ &
$q$ & $\overline{t}_{2}$ &
$q$ & $\overline{t}_{2}^{ \vphantom{H^{L}}}$ \\
 \hline
 61981 & 1407  &  61987 & 1407  &  61991 & 1412  &  62003 & 1404  &  62011 & 1412  &  62017 & 1408  &  62039 & 1407  \\
 62047 & 1408  &  62053 & 1412  &  62057 & 1406  &  62071 & 1407  &  62081 & 1407  &  62099 & 1411  &  62119 & 1409  \\
 62129 & 1412  &  62131 & 1412  &  62137 & 1406  &  62141 & 1413  &  62143 & 1408  &  62171 & 1420  &  62189 & 1409  \\
 62191 & 1416  &  62201 & 1410  &  62207 & 1414  &  62213 & 1414  &  62219 & 1411  &  62233 & 1413  &  62273 & 1413  \\
 62297 & 1412  &  62299 & 1412  &  62303 & 1414  &  62311 & 1411  &  62323 & 1415  &  62327 & 1417  &  62347 & 1411  \\
 62351 & 1412  &  62383 & 1416  &  62401 & 1416  &  62417 & 1414  &  62423 & 1412  &  62459 & 1420  &  62467 & 1419  \\
 62473 & 1419  &  62477 & 1413  &  62483 & 1412  &  62497 & 1416  &  62501 & 1414  &  62507 & 1412  &  62533 & 1413  \\
\hline
 \end{tabular}

 }\newpage

\textbf{Table 5.} Continue 19 \medskip\\\small{
 \renewcommand{\arraystretch}{1.00}
 \begin{tabular}
 {@{}c@{\,\,\,\,\,}c@{\,\,}|@{\,\,}c@{\,\,\,\,\,}c@{\,\,}|@{\,\,}c@{\,\,\,\,\,}c@{\,\,}|@{\,\,}c@{\,\,\,\,\,}c@{\,\,}|@{\,\,}c@{\,\,\,\,\,}c@{\,\,}|@{\,\,}c@{\,\,\,\,\,}c@{\,\,}|@{\,\,}c@{\,\,\,\,\,}c@{}}
  \hline
$q$ & $\overline{t}_{2}$ &
$q$ & $\overline{t}_{2}$ &
$q$ & $\overline{t}_{2}$ &
$q$ & $\overline{t}_{2}$ &
$q$ & $\overline{t}_{2}$ &
$q$ & $\overline{t}_{2}$ &
$q$ & $\overline{t}_{2}^{ \vphantom{H^{L}}}$ \\
 \hline
 65141 & 1450  &  65147 & 1457  &  65167 & 1452  &  65171 & 1452  &  65173 & 1452  &  65179 & 1455  &  65183 & 1449  \\
 65203 & 1450  &  65213 & 1453  &  65239 & 1454  &  65257 & 1452  &  65267 & 1449  &  65269 & 1453  &  65287 & 1458  \\
 65293 & 1459  &  65309 & 1451  &  65323 & 1453  &  65327 & 1448  &  65353 & 1450  &  65357 & 1456  &  65371 & 1454  \\
 65381 & 1457  &  65393 & 1448  &  65407 & 1456  &  65413 & 1454  &  65419 & 1458  &  65423 & 1454  &  65437 & 1455  \\
 65447 & 1451  &  65449 & 1453  &  65479 & 1461  &  65497 & 1460  &  65519 & 1452  &  65521 & 1455  &  65536 & 1453  \\
 65537 & 1457  &  65539 & 1457  &  65543 & 1451  &  65551 & 1459  &  65557 & 1454  &  65563 & 1457  &  65579 & 1454  \\
 65581 & 1459  &  65587 & 1455  &  65599 & 1461  &  65609 & 1455  &  65617 & 1455  &  65629 & 1454  &  65633 & 1463  \\
\hline
 \end{tabular}

 }\newpage

\textbf{Table 5.} Continue 20 \medskip\\\small{
 \renewcommand{\arraystretch}{1.00}
 \begin{tabular}
 {@{}c@{\,\,\,\,\,}c@{\,\,}|@{\,\,}c@{\,\,\,\,\,}c@{\,\,}|@{\,\,}c@{\,\,\,\,\,}c@{\,\,}|@{\,\,}c@{\,\,\,\,\,}c@{\,\,}|@{\,\,}c@{\,\,\,\,\,}c@{\,\,}|@{\,\,}c@{\,\,\,\,\,}c@{\,\,}|@{\,\,}c@{\,\,\,\,\,}c@{}}
  \hline
$q$ & $\overline{t}_{2}$ &
$q$ & $\overline{t}_{2}$ &
$q$ & $\overline{t}_{2}$ &
$q$ & $\overline{t}_{2}$ &
$q$ & $\overline{t}_{2}$ &
$q$ & $\overline{t}_{2}$ &
$q$ & $\overline{t}_{2}^{ \vphantom{H^{L}}}$ \\
 \hline
 68171 & 1492  &  68207 & 1486  &  68209 & 1488  &  68213 & 1487  &  68219 & 1486  &  68227 & 1487  &  68239 & 1489  \\
 68261 & 1489  &  68279 & 1493  &  68281 & 1492  &  68311 & 1492  &  68329 & 1488  &  68351 & 1490  &  68371 & 1494  \\
 68389 & 1497  &  68399 & 1489  &  68437 & 1490  &  68443 & 1485  &  68447 & 1495  &  68449 & 1492  &  68473 & 1489  \\
 68477 & 1493  &  68483 & 1488  &  68489 & 1499  &  68491 & 1487  &  68501 & 1490  &  68507 & 1496  &  68521 & 1493  \\
 68531 & 1493  &  68539 & 1490  &  68543 & 1493  &  68567 & 1495  &  68581 & 1492  &  68597 & 1493  &  68611 & 1498  \\
 68633 & 1495  &  68639 & 1495  &  68659 & 1492  &  68669 & 1490  &  68683 & 1496  &  68687 & 1494  &  68699 & 1492  \\
 68711 & 1495  &  68713 & 1501  &  68729 & 1497  &  68737 & 1493  &  68743 & 1493  &  68749 & 1496  &  68767 & 1493  \\
\hline
 \end{tabular}

 }\newpage

\textbf{Table 5.} Continue 21 \medskip\\\small{
 \renewcommand{\arraystretch}{1.00}
 \begin{tabular}
 {@{}c@{\,\,\,\,\,}c@{\,\,}|@{\,\,}c@{\,\,\,\,\,}c@{\,\,}|@{\,\,}c@{\,\,\,\,\,}c@{\,\,}|@{\,\,}c@{\,\,\,\,\,}c@{\,\,}|@{\,\,}c@{\,\,\,\,\,}c@{\,\,}|@{\,\,}c@{\,\,\,\,\,}c@{\,\,}|@{\,\,}c@{\,\,\,\,\,}c@{}}
  \hline
$q$ & $\overline{t}_{2}$ &
$q$ & $\overline{t}_{2}$ &
$q$ & $\overline{t}_{2}$ &
$q$ & $\overline{t}_{2}$ &
$q$ & $\overline{t}_{2}$ &
$q$ & $\overline{t}_{2}$ &
$q$ & $\overline{t}_{2}^{ \vphantom{H^{L}}}$ \\
 \hline
 71347 & 1524  &  71353 & 1524  &  71359 & 1528  &  71363 & 1522  &  71387 & 1527  &  71389 & 1527  &  71399 & 1531  \\
 71411 & 1527  &  71413 & 1531  &  71419 & 1524  &  71429 & 1532  &  71437 & 1524  &  71443 & 1535  &  71453 & 1529  \\
 71471 & 1526  &  71473 & 1526  &  71479 & 1529  &  71483 & 1533  &  71503 & 1525  &  71527 & 1527  &  71537 & 1527  \\
 71549 & 1528  &  71551 & 1531  &  71563 & 1525  &  71569 & 1528  &  71593 & 1528  &  71597 & 1535  &  71633 & 1523  \\
 71647 & 1532  &  71663 & 1531  &  71671 & 1526  &  71693 & 1535  &  71699 & 1532  &  71707 & 1535  &  71711 & 1531  \\
 71713 & 1527  &  71719 & 1530  &  71741 & 1535  &  71761 & 1531  &  71777 & 1533  &  71789 & 1531  &  71807 & 1528  \\
 71809 & 1530  &  71821 & 1531  &  71837 & 1533  &  71843 & 1532  &  71849 & 1534  &  71861 & 1535  &  71867 & 1535  \\
\hline
 \end{tabular}

 }\newpage

\textbf{Table 5.} Continue 22 \medskip\\\small{
 \renewcommand{\arraystretch}{1.00}
 \begin{tabular}
 {@{}c@{\,\,\,\,\,}c@{\,\,}|@{\,\,}c@{\,\,\,\,\,}c@{\,\,}|@{\,\,}c@{\,\,\,\,\,}c@{\,\,}|@{\,\,}c@{\,\,\,\,\,}c@{\,\,}|@{\,\,}c@{\,\,\,\,\,}c@{\,\,}|@{\,\,}c@{\,\,\,\,\,}c@{\,\,}|@{\,\,}c@{\,\,\,\,\,}c@{}}
  \hline
$q$ & $\overline{t}_{2}$ &
$q$ & $\overline{t}_{2}$ &
$q$ & $\overline{t}_{2}$ &
$q$ & $\overline{t}_{2}$ &
$q$ & $\overline{t}_{2}$ &
$q$ & $\overline{t}_{2}$ &
$q$ & $\overline{t}_{2}^{ \vphantom{H^{L}}}$ \\
 \hline
 74441 & 1562  &  74449 & 1559  &  74453 & 1563  &  74471 & 1562  &  74489 & 1563  &  74507 & 1563  &  74509 & 1569  \\
 74521 & 1566  &  74527 & 1565  &  74531 & 1564  &  74551 & 1567  &  74561 & 1565  &  74567 & 1568  &  74573 & 1561  \\
 74587 & 1570  &  74597 & 1567  &  74609 & 1564  &  74611 & 1562  &  74623 & 1567  &  74653 & 1571  &  74687 & 1568  \\
 74699 & 1571  &  74707 & 1567  &  74713 & 1567  &  74717 & 1563  &  74719 & 1563  &  74729 & 1567  &  74731 & 1573  \\
 74747 & 1568  &  74759 & 1568  &  74761 & 1568  &  74771 & 1568  &  74779 & 1563  &  74797 & 1565  &  74821 & 1563  \\
 74827 & 1567  &  74831 & 1560  &  74843 & 1569  &  74857 & 1567  &  74861 & 1566  &  74869 & 1566  &  74873 & 1571  \\
 74887 & 1574  &  74891 & 1569  &  74897 & 1566  &  74903 & 1570  &  74923 & 1566  &  74929 & 1566  &  74933 & 1570  \\
\hline
 \end{tabular}

 }\newpage

\textbf{Table 5.} Continue 23 \medskip\\\small{
 \renewcommand{\arraystretch}{1.00}
 \begin{tabular}
 {@{}c@{\,\,\,\,\,}c@{\,\,}|@{\,\,}c@{\,\,\,\,\,}c@{\,\,}|@{\,\,}c@{\,\,\,\,\,}c@{\,\,}|@{\,\,}c@{\,\,\,\,\,}c@{\,\,}|@{\,\,}c@{\,\,\,\,\,}c@{\,\,}|@{\,\,}c@{\,\,\,\,\,}c@{\,\,}|@{\,\,}c@{\,\,\,\,\,}c@{}}
  \hline
$q$ & $\overline{t}_{2}$ &
$q$ & $\overline{t}_{2}$ &
$q$ & $\overline{t}_{2}$ &
$q$ & $\overline{t}_{2}$ &
$q$ & $\overline{t}_{2}$ &
$q$ & $\overline{t}_{2}$ &
$q$ & $\overline{t}_{2}^{ \vphantom{H^{L}}}$ \\
 \hline
 77563 & 1604  &  77569 & 1605  &  77573 & 1603  &  77587 & 1597  &  77591 & 1605  &  77611 & 1589  &  77617 & 1595  \\
 77621 & 1602  &  77641 & 1605  &  77647 & 1603  &  77659 & 1603  &  77681 & 1598  &  77687 & 1601  &  77689 & 1605  \\
 77699 & 1597  &  77711 & 1601  &  77713 & 1602  &  77719 & 1601  &  77723 & 1599  &  77731 & 1600  &  77743 & 1603  \\
 77747 & 1604  &  77761 & 1600  &  77773 & 1604  &  77783 & 1604  &  77797 & 1603  &  77801 & 1598  &  77813 & 1605  \\
 77839 & 1602  &  77849 & 1595  &  77863 & 1605  &  77867 & 1607  &  77893 & 1608  &  77899 & 1604  &  77929 & 1609  \\
 77933 & 1608  &  77951 & 1603  &  77969 & 1610  &  77977 & 1606  &  77983 & 1605  &  77999 & 1604  &  78007 & 1602  \\
 78017 & 1603  &  78031 & 1604  &  78041 & 1604  &  78049 & 1602  &  78059 & 1608  &  78079 & 1602  &  78101 & 1603  \\
\hline
 \end{tabular}

 }\newpage

\textbf{Table 5.} Continue 24 \medskip\\\small{
 \renewcommand{\arraystretch}{1.00}
 \begin{tabular}
 {@{}c@{\,\,\,\,\,}c@{\,\,}|@{\,\,}c@{\,\,\,\,\,}c@{\,\,}|@{\,\,}c@{\,\,\,\,\,}c@{\,\,}|@{\,\,}c@{\,\,\,\,\,}c@{\,\,}|@{\,\,}c@{\,\,\,\,\,}c@{\,\,}|@{\,\,}c@{\,\,\,\,\,}c@{\,\,}|@{\,\,}c@{\,\,\,\,\,}c@{}}
  \hline
$q$ & $\overline{t}_{2}$ &
$q$ & $\overline{t}_{2}$ &
$q$ & $\overline{t}_{2}$ &
$q$ & $\overline{t}_{2}$ &
$q$ & $\overline{t}_{2}$ &
$q$ & $\overline{t}_{2}$ &
$q$ & $\overline{t}_{2}^{ \vphantom{H^{L}}}$ \\
 \hline
 80687 & 1636  &  80701 & 1637  &  80713 & 1637  &  80737 & 1638  &  80747 & 1641  &  80749 & 1634  &  80761 & 1632  \\
 80777 & 1637  &  80779 & 1633  &  80783 & 1633  &  80789 & 1635  &  80803 & 1632  &  80809 & 1634  &  80819 & 1636  \\
 80831 & 1634  &  80833 & 1636  &  80849 & 1636  &  80863 & 1637  &  80897 & 1645  &  80909 & 1634  &  80911 & 1637  \\
 80917 & 1638  &  80923 & 1641  &  80929 & 1639  &  80933 & 1635  &  80953 & 1636  &  80963 & 1646  &  80989 & 1638  \\
 81001 & 1642  &  81013 & 1640  &  81017 & 1640  &  81019 & 1637  &  81023 & 1639  &  81031 & 1634  &  81041 & 1641  \\
 81043 & 1637  &  81047 & 1641  &  81049 & 1647  &  81071 & 1638  &  81077 & 1635  &  81083 & 1643  &  81097 & 1631  \\
 81101 & 1635  &  81119 & 1643  &  81131 & 1639  &  81157 & 1639  &  81163 & 1636  &  81173 & 1638  &  81181 & 1641  \\
\hline
 \end{tabular}

 }\newpage

\textbf{Table 5.} Continue 25 \medskip\\\small{
 \renewcommand{\arraystretch}{1.00}
 \begin{tabular}
 {@{}c@{\,\,\,\,\,}c@{\,\,}|@{\,\,}c@{\,\,\,\,\,}c@{\,\,}|@{\,\,}c@{\,\,\,\,\,}c@{\,\,}|@{\,\,}c@{\,\,\,\,\,}c@{\,\,}|@{\,\,}c@{\,\,\,\,\,}c@{\,\,}|@{\,\,}c@{\,\,\,\,\,}c@{\,\,}|@{\,\,}c@{\,\,\,\,\,}c@{}}
  \hline
$q$ & $\overline{t}_{2}$ &
$q$ & $\overline{t}_{2}$ &
$q$ & $\overline{t}_{2}$ &
$q$ & $\overline{t}_{2}$ &
$q$ & $\overline{t}_{2}$ &
$q$ & $\overline{t}_{2}$ &
$q$ & $\overline{t}_{2}^{ \vphantom{H^{L}}}$ \\
 \hline
 83833 & 1678  &  83843 & 1672  &  83857 & 1665  &  83869 & 1673  &  83873 & 1665  &  83891 & 1672  &  83903 & 1675  \\
 83911 & 1666  &  83921 & 1669  &  83933 & 1670  &  83939 & 1672  &  83969 & 1673  &  83983 & 1678  &  83987 & 1675  \\
 84011 & 1678  &  84017 & 1673  &  84047 & 1679  &  84053 & 1671  &  84059 & 1673  &  84061 & 1676  &  84067 & 1669  \\
 84089 & 1672  &  84121 & 1675  &  84127 & 1671  &  84131 & 1674  &  84137 & 1678  &  84143 & 1676  &  84163 & 1671  \\
 84179 & 1676  &  84181 & 1675  &  84191 & 1676  &  84199 & 1678  &  84211 & 1670  &  84221 & 1676  &  84223 & 1677  \\
 84229 & 1671  &  84239 & 1680  &  84247 & 1676  &  84263 & 1679  &  84299 & 1675  &  84307 & 1671  &  84313 & 1671  \\
 84317 & 1679  &  84319 & 1680  &  84347 & 1672  &  84349 & 1673  &  84377 & 1676  &  84389 & 1679  &  84391 & 1674  \\
\hline
 \end{tabular}

 }\newpage

\textbf{Table 5.} Continue 26 \medskip\\\small{
 \renewcommand{\arraystretch}{0.98}
 \begin{tabular}
 {@{}c@{\,\,\,\,\,}c@{\,\,}|@{\,\,}c@{\,\,\,\,\,}c@{\,\,}|@{\,\,}c@{\,\,\,\,\,}c@{\,\,}|@{\,\,}c@{\,\,\,\,\,}c@{\,\,}|@{\,\,}c@{\,\,\,\,\,}c@{\,\,}|@{\,\,}c@{\,\,\,\,\,}c@{\,\,}|@{\,\,}c@{\,\,\,\,\,}c@{}}
  \hline
$q$ & $\overline{t}_{2}$ &
$q$ & $\overline{t}_{2}$ &
$q$ & $\overline{t}_{2}$ &
$q$ & $\overline{t}_{2}$ &
$q$ & $\overline{t}_{2}$ &
$q$ & $\overline{t}_{2}$ &
$q$ & $\overline{t}_{2}^{ \vphantom{H^{L}}}$ \\
 \hline
 87071 & 1700  &  87083 & 1713  &  87103 & 1712  &  87107 & 1706  &  87119 & 1703  &  87121 & 1703  &  87133 & 1706  \\
 87149 & 1708  &  87151 & 1704  &  87179 & 1707  &  87181 & 1707  &  87187 & 1708  &  87211 & 1713  &  87221 & 1701  \\
 87223 & 1703  &  87251 & 1715  &  87253 & 1712  &  87257 & 1709  &  87277 & 1715  &  87281 & 1707  &  87293 & 1710  \\
 87299 & 1703  &  87313 & 1708  &  87317 & 1713  &  87323 & 1715  &  87337 & 1708  &  87359 & 1708  &  87383 & 1712  \\
 87403 & 1708  &  87407 & 1707  &  87421 & 1716  &  87427 & 1711  &  87433 & 1705  &  87443 & 1704  &  87473 & 1701  \\
 87481 & 1703  &  87491 & 1712  &  87509 & 1709  &  87511 & 1711  &  87517 & 1711  &  87523 & 1712  &  87539 & 1706  \\
 87541 & 1714  &  87547 & 1708  &  87553 & 1707  &  87557 & 1712  &  87559 & 1709  &  87583 & 1707  &  87587 & 1713  \\
\hline
 \end{tabular}

 }\newpage

\textbf{Table 5.} Continue 27 \medskip\\\small{
 \renewcommand{\arraystretch}{0.98}
 \begin{tabular}
 {@{}c@{\,\,\,\,\,}c@{\,\,}|@{\,\,}c@{\,\,\,\,\,}c@{\,\,}|@{\,\,}c@{\,\,\,\,\,}c@{\,\,}|@{\,\,}c@{\,\,\,\,\,}c@{\,\,}|@{\,\,}c@{\,\,\,\,\,}c@{\,\,}|@{\,\,}c@{\,\,\,\,\,}c@{\,\,}|@{\,\,}c@{\,\,\,\,\,}c@{}}
  \hline
$q$ & $\overline{t}_{2}$ &
$q$ & $\overline{t}_{2}$ &
$q$ & $\overline{t}_{2}$ &
$q$ & $\overline{t}_{2}$ &
$q$ & $\overline{t}_{2}$ &
$q$ & $\overline{t}_{2}$ &
$q$ & $\overline{t}_{2}^{ \vphantom{H^{L}}}$ \\
 \hline
 90271 & 1740  &  90281 & 1742  &  90289 & 1730  &  90313 & 1735  &  90353 & 1738  &  90359 & 1740  &  90371 & 1736  \\
 90373 & 1740  &  90379 & 1741  &  90397 & 1738  &  90401 & 1745  &  90403 & 1737  &  90407 & 1742  &  90437 & 1747  \\
 90439 & 1741  &  90469 & 1746  &  90473 & 1746  &  90481 & 1736  &  90499 & 1749  &  90511 & 1742  &  90523 & 1741  \\
 90527 & 1743  &  90529 & 1742  &  90533 & 1740  &  90547 & 1744  &  90583 & 1738  &  90599 & 1737  &  90617 & 1731  \\
 90619 & 1747  &  90631 & 1739  &  90641 & 1746  &  90647 & 1740  &  90659 & 1742  &  90677 & 1747  &  90679 & 1745  \\
 90697 & 1738  &  90703 & 1735  &  90709 & 1746  &  90731 & 1740  &  90749 & 1745  &  90787 & 1750  &  90793 & 1745  \\
 90803 & 1746  &  90821 & 1750  &  90823 & 1753  &  90833 & 1746  &  90841 & 1745  &  90847 & 1746  &  90863 & 1748  \\
\hline
 \end{tabular}

 }\newpage

\textbf{Table 5.} Continue 28 \medskip\\\small{
 \renewcommand{\arraystretch}{0.98}
 \begin{tabular}
 {@{}c@{\,\,\,\,\,}c@{\,\,}|@{\,\,}c@{\,\,\,\,\,}c@{\,\,}|@{\,\,}c@{\,\,\,\,\,}c@{\,\,}|@{\,\,}c@{\,\,\,\,\,}c@{\,\,}|@{\,\,}c@{\,\,\,\,\,}c@{\,\,}|@{\,\,}c@{\,\,\,\,\,}c@{\,\,}|@{\,\,}c@{\,\,\,\,\,}c@{}}
  \hline
$q$ & $\overline{t}_{2}$ &
$q$ & $\overline{t}_{2}$ &
$q$ & $\overline{t}_{2}$ &
$q$ & $\overline{t}_{2}$ &
$q$ & $\overline{t}_{2}$ &
$q$ & $\overline{t}_{2}$ &
$q$ & $\overline{t}_{2}^{ \vphantom{H^{L}}}$ \\
 \hline
 93497 & 1772  &  93503 & 1777  &  93523 & 1778  &  93529 & 1779  &  93553 & 1775  &  93557 & 1774  &  93559 & 1778  \\
 93563 & 1770  &  93581 & 1775  &  93601 & 1778  &  93607 & 1781  &  93629 & 1775  &  93637 & 1777  &  93683 & 1779  \\
 93701 & 1772  &  93703 & 1775  &  93719 & 1776  &  93739 & 1778  &  93761 & 1775  &  93763 & 1782  &  93787 & 1776  \\
 93809 & 1780  &  93811 & 1775  &  93827 & 1771  &  93851 & 1771  &  93871 & 1775  &  93887 & 1775  &  93889 & 1776  \\
 93893 & 1777  &  93901 & 1779  &  93911 & 1777  &  93913 & 1777  &  93923 & 1778  &  93937 & 1779  &  93941 & 1778  \\
 93949 & 1781  &  93967 & 1781  &  93971 & 1786  &  93979 & 1778  &  93983 & 1773  &  93997 & 1779  &  94007 & 1779  \\
 94009 & 1775  &  94033 & 1774  &  94049 & 1785  &  94057 & 1780  &  94063 & 1784  &  94079 & 1778  &  94099 & 1784  \\
\hline
 \end{tabular}%

 }\newpage

\textbf{Table 6.} The smallest known sizes
 $\overline{t}_{2}=\overline{t}_{2}(2,q)$
 of complete arcs in planes
 $\mathrm{PG}(2,q),$\\ $100001\leq q\leq 160000$, $q$ power prime\footnote{For
 $q\ge150011$, all the sizes are new in comparison with
 \cite{BDFKMP-ArXiv}. They are obtained in \cite{BDFKMP-JGtoappear}.}\medskip\\\small{
 \renewcommand{\arraystretch}{0.92}
 %

 \end{document}